\documentclass{article}
\usepackage[utf8]{inputenc}

\usepackage{amsmath,amssymb,amsthm}
\usepackage{cleveref,color}

\usepackage{graphicx}

\usepackage[sort,compress,square,numbers]{natbib}
\crefname{equation}{}{}
\crefname{figure}{Figure}{Figures}
\crefname{section}{Section}{Sections}
\crefname{table}{Table}{Tables}
\crefname{appendix}{Appendix}{Appendices}

\DeclareMathOperator*{\argmin}{arg\,min}
\DeclareMathOperator*{\locargmin}{loc\,arg\,min}

\DeclareMathOperator{\imag}{Im}
\DeclareMathOperator{\real}{Re}

\newcommand{\bbR}{\mathbb{R}}
\newcommand{\bbRnn}{\mathbb{R}_{\geq 0}}
\newcommand{\bbC}{\mathbb{C}}
\newcommand{\uscat}{u^{\textrm{scat}}}
\newcommand{\umeas}{{\bf u}^{\textrm{meas}}}
\newcommand{\uinc}{u^{\textrm{inc}}}
\newcommand{\uout}{{\bf u}^{\textrm{out}}}
\newcommand{\xb}{{\mathbf x}}
\newcommand{\wb}{{\mathbf w}}
\newcommand{\zb}{{\mathbf z}}
\newcommand{\fb}{{\mathbf f}}
\newcommand{\Jb}{{\mathbf J}}

\newcommand{\nb}{{\mathbf n}}
\newcommand{\hb}{{\mathbf h}}
\newcommand{\yb}{{\mathbf y}}
\newcommand{\rb}{{\mathbf r}}
\newcommand{\db}{{\mathbf d}}
\newcommand{\cb}{{\mathbf c}}
\newcommand{\vb}{{\mathbf v}}
\newcommand{\eb}{{\mathbf e}}
\newcommand{\gammab}{{\boldsymbol \gamma}}
\newcommand{\alphab}{{\boldsymbol \alpha}}
\newcommand{\betab}{{\boldsymbol \beta}}
\newcommand{\etab}{{\boldsymbol \eta}}

\newcommand{\Ccal}{\mathcal{C}}
\newcommand{\cK}{\mathcal{K}}
\newcommand{\cS}{\mathcal{S}}
\newcommand{\cD}{\mathcal{D}}
\newcommand{\cT}{\mathcal{T}}

\newcommand{\Ccalch}{\mathcal{C}_{\textrm{CH}}}
\newcommand{\Ccalabv}{\mathcal{C}_{\textrm{ABV}}}
\newcommand{\zerob}{{\mathbf 0}}

\newcommand{\ktwomax}{k_2^{\textrm{max}}}

\newcommand{\Fimp}{\mathcal{F}^{\textrm{imp}}}
\newcommand{\Fneu}{\mathcal{F}^{\textrm{neu}}}
\newcommand{\Ftrans}{\mathcal{F}^{\textrm{trans}}}

\newcommand{\ccurv}{C_{H}}

\newcommand{\lamfs}{{\lambda_{\textrm{FS}}}}
\newcommand{\lamch}{{\lambda_{\textrm{CH}}}}
\newcommand{\lamabv}{{\lambda_{\textrm{ABV}}}}

\newcommand{\gn}{ {\textrm{gn}}}
\newcommand{\sd}{ {\textrm{sd}}}
\newcommand{\slf}{ {\textrm{sf}}}
\newcommand{\gf}{ {\textrm{gf}}}
\newcommand{\filt}{ {\textrm{filt}}}

\DeclareMathOperator{\proj}{proj}
\DeclareMathOperator{\area}{area}
\DeclareMathOperator{\fp}{f.p.}

\renewcommand{\bar}{\overline}

\newcommand{\eps}{\epsilon}

\newtheorem{theorem}{Theorem}
\newtheorem{proposition}{Proposition}
\newtheorem{remark}{Remark}

\title{
Reconstructing the shape and material parameters of dissipative obstacles using an impedance model  }
\author{Travis Askham\thanks{Department of Mathematical Sciences, New Jersey Institute of Technology, Newark, NJ, USA. \textit{Email: askham@njit.edu}}
\and
Carlos Borges\thanks{Department of Mathematics, University of Central Florida, Orlando, FL, USA. \textit{Email: carlos.borges@ucf.edu}}}

\begin{document}

\maketitle

\begin{abstract}
  In inverse scattering problems, a model that allows for the simultaneous
  recovery of both the domain shape and an impedance boundary condition covers
  a wide range of problems with impenetrable domains, including recovering the
  shape of sound-hard and sound-soft obstacles and obstacles with thin coatings.
  This work develops an optimization framework for recovering the shape and material parameters
  of a penetrable, dissipative obstacle in the multifrequency setting, using a
  constrained class of curvature-dependent impedance function models proposed by Antoine,
  Barucq, and Vernhet~\cite{antoine2001high}. We find that
  this constrained model improves the robustness of the recovery problem, compared
  to more general models, and provides meaningfully better obstacle recovery than
  simpler models. We explore the effectiveness of the model for varying
  levels of dissipation, for noise-corrupted data, and for limited aperture
  data in the numerical examples. 
\end{abstract}

\section{Introduction}

Inverse scattering problems arise in many applications, including
sensing \cite{ustinov2014geophysical,jin2006theory}, ocean acoustics
\cite{chavent2012inverse,collins1994inverse}, medical imaging
\cite{kuchment2014radon,nashed2002inverse,scherzer2010handbook,simonetti2008inverse},
nondestructive testing \cite{collins1995nondestructive,engl2012inverse,langenberg1993imaging},
and radar and sonar \cite{cheney2009fundamentals}. The general setting of
those problems is characterized by probing a domain with one or multiple incident
waves and obtaining measurements of the scattered waves. From this scattered wave data,
one seeks to reconstruct some property of the domain, e.g. density, sound-speed,
impedance, shape, etc.

In this work, we consider the problem of recovering the shape and physical
parameters of homogeneous, penetrable
obstacles from measurements of the scattered field in the far field region through
the use of an impedance model approximation of the standard transmission model. We consider
dissipative obstacles in which acoustic waves in the medium are damped and acoustic energy
dissipates into thermal energy~\cite{blauert2009dissipation}. In the frequency domain, the
dissipation corresponds to a complex wavenumber for the medium of the form
$k_1=\omega\sqrt{(1+\imath \delta/\omega)}/c_1$, where $\omega$ is the pulsating frequency of the
incident wave, $c_1$ is the sound speed of the medium, and $\delta>0$ is the dissipation
constant for the medium.

Let $\Omega_1$ denote the interior of a dissipative obstacle with boundary
$\Gamma$ and let ${\bf n}$ denote the outward normal on $\Gamma$. Given the incident field, $\uinc$, the scattered field, $\uscat$, is modeled by
the Helmholtz equation with a transmission boundary condition. In the notation
of~\cite{antoine2005construction}, we have
\begin{equation}\label{eq:transpde}
\begin{aligned}
    -(\Delta + k_2^2)\uscat &= 0 & \textrm{ in } \Omega_2 \; , \\
    -(\Delta + k_1^2)\uscat &= k_2^2(1-N^2) \uinc & \textrm{ in } \Omega_1 \; , \\
    \left [ \uscat \right ] &= 0 & \textrm{ on } \Gamma \; , \\
    \left [ \chi \partial_n \uscat \right ] &=
    -\left [ \chi \partial_n \uinc \right ] & \textrm{ on } \Gamma \; , \\
    \sqrt{|\xb|} \left(\uscat - \imath k_2\frac{\xb}{|\xb|}\cdot
    \nabla \uscat \right) &\to 0 & \textrm{ as } |\xb|\to \infty \; , 
\end{aligned}
\end{equation}
where $\Omega_2 = \bbR^2 \setminus \bar{\Omega_1}$; $k_2=\omega/c_2$ is the wavenumber of the
incoming incident wave; $c_2$ is the sound speed of the background medium $\Omega_2$
and $c_r=c_1/c_2$ is the relative sound speed;
$N=\sqrt{(1+\imath\delta/\omega)}/c_r$ is the relative refractive index;
$\rho_1$ and $\rho_2$ are the densities for $\Omega_1$ and $\Omega_2$, respectively,
and $\rho_r = \rho_1/\rho_2$;
$\alpha=1/(\rho_r(1+\imath\delta/\omega))$ is the complex contrast coefficient;
and the function $\chi$ is equal to $1$ in $\Omega_2$ and $\alpha$ in $\Omega_1$.
The notation $\left [\phi \right]$ denotes the difference between the exterior
and interior traces, or the ``jump'', of the function $\phi$ across $\Gamma$.

The transmission problem \cref{eq:transpde}
can be solved using standard boundary integral equation methods (BIEMs); see
\cref{sec:bieform-trans} for details. Like the case of impenetrable obstacles,
say with sound-soft or sound-hard boundary conditions, the solution of a penetrable
transmission problem can then be discretized using unknowns on the boundary of
the domain alone. However, the transmission problem generally requires the solution
of a system of twice the size.

\begin{remark}
  The difference in computational effort between the penetrable and
  impenetrable cases is even more significant in the inverse obstacle
  setting. The Fr\'{e}chet derivative of the PDE solution with respect to
  the material parameters typically involves the solution of an
  inhomogeneous PDE for the penetrable case; see, e.g., \cite[\S A.2]{carpio2020bayesian}.
\end{remark}

In an attempt to decrease the computational and memory costs in the solution of
the problem, several authors proposed to approximate the forward transmission
problem by a forward scattering problem for an impenetrable obstacle with a
generalized impedance boundary condition (GIBC), see
\cite{rytov1940calcul, senior1960impedance, senior1995approximate, jones1992improved, senior1995generalized, senior1995approximateelectro, senior1997higher, wang1987limits, antoine2005approximation, antoine2001high, haddar2004asymptotic}.
Such approximations hold in various settings, including scattering from
penetrable, dissipative objects and objects with thin, penetrable coatings.

In \cite{antoine2001high}, the authors observe that for large dissipation
constant, i.e. $\delta \gg 1$, the wave does not penetrate the medium
significantly and they show that the forward transmission problem can be approximated
asymptotically by a related problem with a GIBC. In particular, if $\uscat$
is the solution of
\cref{eq:transpde}, then
there exists a local operator $\mathcal{Y}$ such that $\uscat \approx \phi$, where
\begin{equation} \label{eq:imppde}
\begin{aligned}
  -(\Delta + k_2^2)\phi &= 0 & \textrm{ in } \Omega_2 \; , \\
  \left ( \partial_n - \mathcal{Y} \right )\phi &=
  -\left ( \partial_n - \mathcal{Y} \right ) \uinc \; & \textrm{ on } \Gamma \; , \\
  \sqrt{|\xb|} \left(\phi - ik_2\frac{\xb}{|\xb|}\cdot
  \nabla \phi \right) &\to 0 & \textrm{ as } |\xb|\to \infty \; .
\end{aligned}
\end{equation}

\begin{figure}[ht]
  \center
  \includegraphics[width=0.9\textwidth]{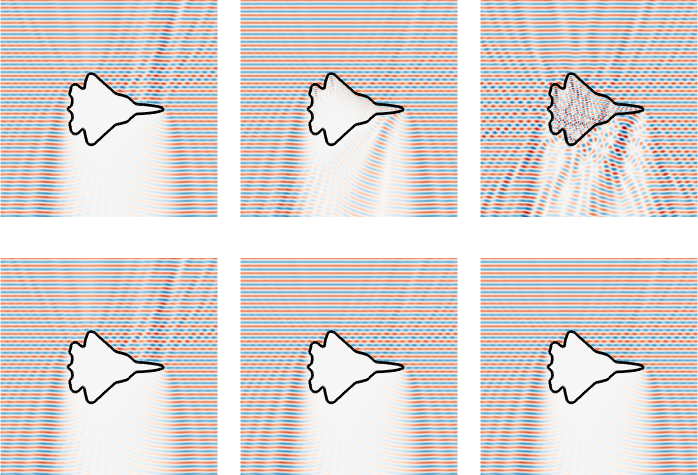}
  \caption{Comparison of the total field 
    of the solution of the transmission ($u=\uscat+\uinc$) and impedance ($u=\phi+\uinc$)
    scattering problems for $\omega=30$. The
    top row shows the transmission solution and the bottom row shows the solution of the
    first order impedance model. The dissipation decreases from left to right.
    The first column sets $\delta = \delta_0 = \sqrt{3} \omega$.
    The second and third columns have $\delta = \delta_0/16$ and $\delta = \delta_0/256$,
    respectively.}
  \label{fig:impvstrans-scatter}
\end{figure}


It is possible to obtain different order approximations for the operator
$\mathcal{Y}$. For example, the first order operator derived in~\cite{antoine2001high}
is a multiplication operator of the form
\begin{equation}
  \label{eq:abvorder1}
  \left [ \mathcal{Y}_1 \phi \right ] (\xb)
  = -\alpha(\imath k_2 N + H(\xb))\phi(\xb) \; ,
\end{equation}
where $H(\xb)$ is the signed curvature at $\xb\in \Gamma$. Higher order GIBCs
involve differentiation operators on the curve~\cite{antoine2001high,haddar2004asymptotic}. 

\Cref{fig:impvstrans-scatter} compares the solution of \cref{eq:transpde}
to the solution of \cref{eq:imppde} with $\mathcal{Y}=\mathcal{Y}_1$
for the same incident field and varying dissipation. With dissipation in the regime
recommended by~\cite{antoine2005construction} ($\delta=\delta_0=\sqrt{3}\omega$), the impedance
model agrees well with the transmission model. Slightly below this regime ($\delta = \delta_0/16$),
the models still agree well in much of the exterior; however, the solution of the
transmission problem exhibits a wave transmitting through the narrowest
part of the plane, which is not captured by the impedance model.
For even lower values of dissipation ($\delta = \delta_0/256$), the solution of
the transmission problem penetrates
the obstacle and the solutions of the two models are notably different.

We will explore the suitability of a simple multiplicative GIBC
of the form $[\mathcal{Y} \phi](\xb) = -\imath k_2 \lambda(\xb;\omega) \phi(\xb)$ as an
approximate forward model in the inverse obstacle scattering setting with
multiple frequency data. More concretely, we consider scattering data for a fixed
set of material parameters: $c_1$, $c_2$, $\rho_1$, $\rho_2$, and $\delta$. We assume that
scattering measurements are available for $N_k$ pulsation frequencies,
$\omega_1 < \omega_2 < \cdots  < \omega_{N_k}$. The incident fields are plane waves
of the form $\uinc(\xb) = \exp(\imath k_2 \db \cdot \xb)$ and data are
available for $N_d(\omega)$ angles of incidence at each $\omega$, with
directions determined by the unit vectors $\db_1,\ldots,\db_{N_d}$.
The scattered fields are measured for each incident wave at $N_r(\omega)$
receptor locations, $\rb_1,\ldots,\rb_{N_r}$, far from the obstacle.

Let the forward transmission operator be defined by
$\Ftrans_{\omega,\db}(\Gamma) = \uout \in \bbC^{N_r}$, where
$\uout_j = \uscat(\rb_j)$ and $\uscat$ is the solution of the transmission
scattering problem, \cref{eq:transpde}.
Similarly, let the forward impedance operator be defined by
$\Fimp_{\omega,\db}(\Gamma,\lambda) = \uout \in \bbC^{N_r}$, where
$\uout_j = \phi(\rb_j)$ and $\phi$ is the solution of the impedance scattering
problem, \cref{eq:imppde}, with $\mathcal{Y}=-\imath k_2\lambda$.

Let $\Gamma_\star$ be the true boundary curve of interest and let
$\umeas_{\omega,\db} = \Ftrans_{\omega,\db}(\Gamma_\star)$. At a given frequency $\omega_j$,
a natural definition for the ``best-fit'' boundary curve and GIBC is then given as the
solution of the following constrained optimization problem:

\begin{equation}
  \label{eq:invprob_orig} 
 \left[\hat{\Gamma}_j,\hat{\lambda}_j \right]= \argmin_{\Gamma \in S_\Gamma(\omega_j), \lambda \in S_\lambda(\omega_j)} \sum_{i=1}^{N_d} | \umeas_{\omega_j,\db_i} -
  \Fimp_{\omega_j,\db_i}(\Gamma,\lambda)|^2 \; ,
\end{equation}
where the sets $S_\Gamma(\omega)$ and $S_\lambda(\omega)$ are chosen to be appropriate
spaces for the curve and GIBC at a given pulsation, and can be designed to
regularize the problem. 

The problems defined by \cref{eq:invprob_orig} are generally nonlinear,
non-convex, and ill-posed.
Following the continuation-in-frequency approach~\cite{chen1997inverse, bao2004inverse, bao2007inverse, borges2017high,sini2012inverse, borges2015inverse},
we solve these problems in sequence beginning at the lowest frequency, where the
problem is approximately convex, and use the solution $\hat{\Gamma}_{j-1},\hat{\lambda}_{j-1}$
as an initial guess for $\hat{\Gamma}_j,\hat{\lambda}_j$. This helps to
mitigate the non-convexity. We apply standard iterative methods to handle
the non-linearity and we select the sets $S_\Gamma$ and $S_\lambda$ to mitigate
the ill-posedness. We then take the recovered boundary and impedance function
to be the solution of \cref{eq:invprob_orig} for the highest frequency data available.

We provide details of the constraint sets,
impedance models, and gradient formulas in \cref{sec:analytical_app}. We
describe the constraint set for the geometry, $S_\Gamma$, in detail in
\cref{sec:cif_constraint_sets}. We present
three options for describing the impedance function $\lambda$ (and
appropriate constraint sets $S_\lambda$) in
\cref{sec:impedance_funs}: a general model based on a Fourier series in arc-length and
two more-constrained models which depend on the curvature function on $\Gamma$.
While Fr\'{e}chet derivative formulae are known for the impedance boundary
value problem, the curvature-dependent models require some new expressions which we
derive in \cref{sec:frechetderivative}.

We describe some details of the iterative
optimization framework we use at each frequency in \cref{sec:optimization_methods}.
The curvature-dependent impedance models are more amenable to the imposition
of physical constraints, which can be handled effectively by projected gradient
methods as described in \cref{sec:projgrad}.

Numerical results are presented in \cref{sec:numerical}. These indicate that
the solution scheme has much greater success with the curvature-dependent models.
Domain recovery is effective well below the level of dissipation needed for
qualitative agreement between the impedance and transmission forward problems
and the curvature-dependent models provide a meaningful advantage over some
simpler alternatives. We find that the dissipation, $\delta$, and the
product $c_r\rho_r$ can be recovered reliably, with sufficient dissipation, but
recovering the values $c_r$ and $\rho_r$ individually appears to be
difficult.

We discuss some implications of these results and possible future directions
in \cref{sec:conclusions}.

\subsection{Relation to the literature}
The inverse problem of recovering the shape and boundary conditions using single frequency
data and a model with the classical impedance boundary conditions was studied by several
authors \cite{ivanyshyn2011inverse, kress2001inverse, serranho2006hybrid, lee2007inverse, akduman2003direct, qin2012inverse, kress2018inverse, smith1985inverse}.
The single frequency inverse scattering problem using the generalized impedance boundary
condition model was considered in \cite{bourgeois2010identification, yaman2019reconstruction, yang2014reconstruction, kress2019some, kress2018integral, cakoni2012integral, bourgeois2012simultaneous, bourgeois2011stable, guo2015multilayered, aslanyurek2014reconstruction}.
The use of multifrequency data by applying continuation in frequency to recover the sound
speed of a volume was studied in \cite{chen1997inverse, bao2004inverse, bao2007inverse, borges2017high},
to recover the shape of an impenetrable domain was studied in \cite{sini2012inverse, borges2015inverse},
and to recover the classical impedance boundary condition simultaneously with the shape of
obstacle was studied in \cite{borges2022multifrequency}. 



The present work introduces a new framework to simultaneously recover the
shape of an obstacle and an appropriate impedance function for multifrequency transmission
data from a dissipative obstacle. We build on techniques previously presented in
\cite{borges2022multifrequency}, in
particular the use of continuation-in-frequency for the inverse problem and
high-order-accurate methods for solving
the associated PDEs via integral equation representations. We provide some necessary formulae
and identify efficient numerical schemes for treating curvature-dependent impedance models, like
the first order model, \cref{eq:abvorder1}, of \cite{antoine2001high}, in the inverse
scattering setting.

\section{Details of the model and gradient formulas}

\label{sec:analytical_app}

This section describes some details of the problem discretization and
the constraint sets which regularize the problem at each frequency. The
discretization of the obstacle boundary curve and appropriate constraints are
described in \cref{sec:cif_constraint_sets}. A couple of competing models
for the impedance function and appropriate constraints for these are described in
\cref{sec:impedance_funs}. Two of these models have the impedance function
depend on the curvature of the obstacle; to apply an iterative method for the
solution of \cref{eq:invprob_orig} with these models, we then require some
new derivative formulas that we derive in \cref{sec:frechetderivative}. 

\subsection{Representation of the obstacle and its constraints}
\label{sec:cif_constraint_sets}

We represent a boundary curve, $\Gamma$, of length $L$ by an arc-length parameterization
$\gammab:[0,L)\to\mathbb{R}^2$, where $\gammab(s)=\left(x(s),y(s)\right)$,
with $x,y:[0,L)\to\mathbb{R}$, being trigonometric polynomials of the form
\begin{equation}
\begin{aligned} \label{eq:domain_coord_func}
  x(s)&=a_{1,0}+\sum_{m=1}^{N(\omega)}\left(a_{1,m}\cos(2\pi ms/L)+b_{1,m}\sin(2\pi ms/L)\right) \;, \\
  y(s)&=a_{2,0}+\sum_{m=1}^{N(\omega)}\left(a_{2,m}\cos(2\pi ms/L)+b_{2,m}\sin(2\pi ms/L)\right) \; ,  
\end{aligned}
\end{equation}
where $a_{j,0}$ and $a_{j,m}$, and $b_{j,m}$ for $j=1,2$ and $m=1,\ldots,N(\omega)$ are real constants
and $N(\omega)$ is an integer proportional to $\omega$.

To ensure that the inverse problem at frequency $\omega$ is well-posed, we
require that the arc-length parameterization of the curve, $\gammab$, have bandlimited
curvature in a suitable sense.
The signed curvature $H$ for a curve-parameterization $\gammab$ is defined as
\begin{equation}
  \label{eq:signedcurve}
H(\gammab) = \frac{x'y'' - x''y'}{(x'^2+y'^2)^{3/2}} \; ,
\end{equation}
where the denominator is constant equal to 1 for an arc-length parameterization. 
The curvature $H(\gammab)$ then has a Fourier expansion of twice the length, i.e. of the form
\begin{equation*}
  H(\gammab)(s) = a_{H,0} + \sum_{m=1}^{2N} \left ( a_{H,m} \cos(2\pi ms/L) +
  b_{H,m} \sin(2\pi ms/L) \right )\; ,
\end{equation*}
where $a_{H,0}$ and $a_{H,m}$ and $b_{H,m}$ for $m=1,\ldots,2N$ are real constants. 
Define $\mathcal{E}(\Gamma)$, and $\mathcal{E}^{M(\omega)}(\Gamma)$ to be the elastic energy
of the curve $\Gamma$ and the elastic energy contained in the first $M(\omega)$ modes
of the curvature, respectively, i.e. 
\begin{equation}
   \label{eq:em}   
\begin{aligned}
  \mathcal{E}(\Gamma) &= \frac{a_{H,0}^2}{2} + \sum_{m=1}^{2N} \left ( a_{H,m}^2 + b_{H,m}^2 \right )  \quad\text{and} \\
  \mathcal{E}^{M(\omega)}(\Gamma) &= \frac{a_{H,0}^2}{2} + \sum_{m=1}^{M(\omega)}
  \left ( a_{H,m}^2 + b_{H,m}^2 \right )\, .
\end{aligned}
\end{equation}
Selecting a value of $M(\omega)$ proportional to $\omega$ and a value of $\ccurv$ near 1,
we can impose a bandlimited curvature requirement as the constraint such that
$\mathcal{E}^{M(\omega)} \geq \ccurv \mathcal{E}(\Gamma)$.

In addition to this constraint on the curvature, we impose that the curve
is simple, i.e. non-self-intersecting. The constraint set for the boundary
curves is then
\begin{equation*}
  \mathcal{S}_\Gamma(\omega) = \{ \Gamma \, |  \, \text{$\Gamma$ is simple
    and } \mathcal{E}^{M(\omega)}(\Gamma) \geq \ccurv \mathcal{E}(\Gamma) \} \, .
\end{equation*}

Beyond the constraint set, the problem can also be regularized by
limiting the possible search directions. Because a tangential update
of the curve does not change the shape, curve updates will always be
in the normal direction, i.e. updates will be
of the form $\gammab(s) \to \gammab(s) + h(s)\nb(s)$ where
$\nb = (y',-x')$ is the normal to the curve.
We will also only propose curve updates with frequency content
proportional to the frequency of the data, i.e. $h:[0,L)\to \bbR$ will
be of the form 
\begin{equation} \label{eq:hupdate}
h(s) = a_{h,0} + \sum_{m=1}^{N_\gamma(\omega)}\left( a_{h,m} \cos(2\pi ms/L)+ b_{h,m} \sin(2\pi ms/L)\right)
\end{equation}
where $a_{h,0}$ and $a_{h,m}$ and $b_{h,m}$ for $m=1,\ldots,N_\gamma(\omega)$ are real
constants and $N_\gamma(\omega)$ is an integer proportional to $\omega$.

\begin{remark}
  It should be noted that the integer parameters in this section,
  $N(\omega)$, $M(\omega)$ and $N_\gamma(\omega)$, are all selected proportional
  to $\omega$, but these values are selected with different goals.
  
  The values of $M(\omega)$ and $N_\gamma(\omega)$ are chosen with respect
  to the scale of boundary features which can stably be recovered given
  scattering data collected at frequency $\omega$. In light of the Heisenberg
  uncertainty principle~\cite{chen1997inverse}, it is unreasonable to
  attempt to stably reconstruct features smaller than half of a wavelength in size.
  
  The value of $N(\omega)$ determines the number of discretization nodes
  which are used to represent the curve and should be selected to be sufficiently
  large that the boundary integral equation method used to approximate the
  solution of the forward problems is accurate. For accuracy, it is typically
  sufficient to sample the curve at some fixed number of points per wavelength.
  One contrast between $N(\omega)$ and the other parameters is that there
  is not much harm in selecting $N(\omega)$ too large, other than the
  unnecessary computational burden.
\end{remark}

\subsection{Impedance function models}
\label{sec:impedance_funs}

We consider three models for the impedance function. The first
is to model the impedance function as a Fourier series in arc-length:

\begin{equation}
  \label{eq:lamfsmodel}
  \lamfs[\cb](s) = \sum_{m=-N_c(\omega)}^{N_c(\omega)} c_m \exp(2\pi \imath m s/L) \; ,
\end{equation}
where $N_c(\omega)$ should be chosen to effectively
regularize the model. The inverse problem
at a single frequency, \cref{eq:invprob_orig}, can then be rephrased as

\begin{equation}
  \label{eq:invprob_fs} 
 \left[\hat{\Gamma}_j,\hat{\cb}_j\right] = \argmin_{\Gamma \in S_\Gamma(\omega_j), \cb \in \bbC^{2N_c(\omega_j)+1}} \sum_{i=1}^{N_d} | \umeas_{\omega_j,\db_i} -
  \Fimp_{\omega_j,\db_i}(\Gamma,\lamfs[\cb])|^2 \; .
\end{equation}

The second approach is to model the impedance as a linear function of
the curvature:
\begin{equation}
  \label{eq:lamchmodel}
  \lamch[\alphab] = \alpha_1 + \alpha_2 H \; ,
\end{equation}
where $\alphab=(\alpha_1,\alpha_2) \in \bbC^2$. The inverse problem at a single frequency,
\cref{eq:invprob_orig}, can then be rephrased as

\begin{equation}
  \label{eq:invprob_ch} 
 \left[\hat{\Gamma}_j,\hat{\alphab}_j\right] = \argmin_{\Gamma \in S_\Gamma(\omega_j), \alphab \in \bbC^{2}} \sum_{i=1}^{N_d} | \umeas_{\omega_j,\db_i} -
  \Fimp_{\omega_j,\db_i}(\Gamma,\lamch[\alphab])|^2 \; .
\end{equation}

We also consider a more restricted form of the curvature-dependent
model which is based directly on the first order model in \cref{eq:abvorder1}:

\begin{equation}
  \label{eq:lamabvmodel}
  \lamabv[\betab] = \beta_2 \sqrt{1-\imath \beta_1}
  - \frac{\imath \beta_3 (1-\imath \beta_1)}{k_2} H \; ,
\end{equation}
where the parameters can be taken to be non-negative and real, i.e.
$\betab \in \bbRnn^3$. The $\betab$ parameters are related to
the physical parameters as follows:

\begin{equation}
\label{eq:betatophys} \beta_1 = \frac{\delta}{\omega} \;, \quad
\beta_2 = \frac{1}{\rho_rc_r\sqrt{1+\delta^2/\omega^2}} \; ,
\quad \beta_3 = \frac{1}{\rho_r(1+\delta^2/\omega^2)} \; .
\end{equation}
The inverse problem at a single frequency \cref{eq:invprob_orig} can
then be rephrased as

\begin{equation}
  \label{eq:invprob_abv} 
 \left[ \hat{\Gamma}_j,\hat{\betab}_j \right]= \argmin_{\Gamma \in S_\Gamma(\omega_j), \betab \in \bbRnn^{3}}
  \sum_{i=1}^{N_d} | \umeas_{\omega_j,\db_i} -
  \Fimp_{\omega_j,\db_i}(\Gamma,\lamabv[\betab])|^2 \; .
\end{equation}

The $\lamabv$ model has
the natural physical constraint that $\betab=(\beta_1,\beta_2,\beta_3) \in \bbRnn^3$. 
Some natural constraints on the $\lamch$ model are that
$\imag(\alpha_1)\leq 0$ and $\imag(\alpha_2)\leq 0$.
For these simple models, it is relatively easy to impose such constraints using
the methods of \cref{sec:projgrad}. Physical constraints for the $\lamfs$
model are both less obvious and more difficult to impose. In the numerical
examples in this manuscript, the $\lamfs$ model coefficients are unconstrained.

\begin{remark}
  Above, we parameterize $\lamabv$ by $\betab$ instead of the
  actual physical parameters, $\delta$, $c_r$, and $\rho_r$. One reason
  for this is that the natural constraints are $\delta \geq 0$,
  $c_r>0$, and $\rho_r > 0$ and the impedance function becomes infinite if
  $c_r = 0$ or $\rho_r = 0$. The constrained optimization problem is
  simpler for the closed, convex constraint set of the $\betab$ variables.
  In particular, we do not have to select arbitrary lower bounds for $c_r$ and
  $\rho_r$.
\end{remark}

\subsection{Fr\'{e}chet derivatives of the forward operator $\Fimp$} \label{sec:frechetderivative}

Let $\Fimp_{\omega,\db}$ be the forward map defined as in the
introduction, i.e. $\Fimp_{\omega,\db}(\Gamma,\lambda) \in \bbC^{N_r}$ is the solution of the
impedance scattering problem \cref{eq:imppde}
for the incident field $\uinc = \exp(\imath k_2\db\cdot \xb)$
evaluated at the receptors, $\rb_j$.
To apply standard iterative solvers to the minimization problem
\cref{eq:invprob_orig}, we require expressions for the derivatives of
$\Fimp_{\omega,\db}$ with respect to the curve $\Gamma$ and the impedance function
$\lambda$. Below, we use the notation $D_{\Gamma}$ to denote the Fr\'{e}chet
derivative of a quantity with respect to the boundary and $D_{\lambda}$ to
denote the Fr\'{e}chet derivative with respect to the impedance function.
We will also drop the dependence of $\Fimp$ on $\omega$ and $\db$ for
the sake of brevity, and, in some cases, even the dependence on $\Gamma$ and $\lambda$
when it is clear.

To be more concrete, let $g(s)$ be a sufficiently smooth complex-valued
function of arc-length on $\Gamma$ and let $\lambda_\eps(s) := \lambda(s)
+ \eps g(s)$ for small $\eps > 0$ be a perturbation of $\lambda$ in the direction $g$.
The Fr\'{e}chet derivative $D_\lambda\Fimp$ is the linear operator such that 

\begin{equation*}
  \lim_{\eps \to 0}\frac{\|\Fimp(\Gamma,\lambda_\eps)
    -\left(\Fimp(\Gamma,\lambda)+ \left [D_\lambda\Fimp \right ] g \right)\|}{\eps} = 0
\end{equation*}
for all such $g$.

Likewise, letting $\hb(s)$ be a sufficiently smooth $\bbR^2$-valued function
of arc-length on $\Gamma$
and $\gammab(s)$ be an arc-length parameterization of $\Gamma$,
we can define a new curve parameterization (no longer in arc-length) by
$\gammab_\eps(s) = \gammab(s) + \eps \hb(s)$ and denote the corresponding
curve by $\Gamma_\eps$. The Fr\'{e}chet derivative $D_\Gamma \Fimp$ is the
linear operator such that 

\begin{equation*}
  \lim_{\eps \to 0}\frac{\|\Fimp(\Gamma_\eps,\lambda)-
    \left(\Fimp(\Gamma,\lambda)+ \left [D_\Gamma\Fimp \right ] \hb\right)\|}{\eps} = 0
\end{equation*}
for all such $\hb$.


Explicit formulas for these derivatives are known. They are expressed
as the solutions of impedance scattering problems:
\begin{theorem}[\cite{hettlich1998frechet,hettlich1995frechet}]
Let $u=\uscat+\uinc$ be the total field for the solution of the original
obstacle problem \cref{eq:imppde}.
The $j$th entry of $\left [D_\lambda\Fimp \right ] g$ is equal to $v_\lambda(\rb_j)$, where 
\begin{equation}\label{eq:imp_derimp_pde}
\begin{aligned}
  -(\Delta + k_2^2)v_\lambda &= 0 & \textrm{ in } \Omega \; , \\
    \partial_n v_\lambda +\imath k_2\lambda v_\lambda &= -\imath k_2 g u & \textrm{ on } \Gamma \; , \\
    \sqrt{|x|} \left(v_\lambda - \imath k_2\frac{x}{|x|}\cdot
    \nabla v_\lambda \right) &\to 0 & \textrm{ as } |x|\to \infty \; . 
\end{aligned}
\end{equation}

The $j$th entry of $\left [D_\Gamma\Fimp \right] \hb$ is equal to $v_\gamma(\rb_j)$, where 
\begin{equation} \label{eq:imp_derdom_pde}
\begin{aligned}
    -(\Delta + k_2^2)v_\gamma &= 0 & \textrm{ in } \Omega \; , \\
    \partial_n v_\gamma +\imath k_2\lambda v_\gamma &= 
    k_2^2 \hb\cdot \nb u + \frac{d}{ds}\left( \hb\cdot \nb \frac{du}{ds} \right) \\
    & \qquad -\lambda \hb \cdot \nb \left(\partial_n u - H u\right) & \textrm{ on } \Gamma \; , \\
    \sqrt{|x|} \left(v_\gamma - \imath k_2\frac{x}{|x|}\cdot
    \nabla v_\gamma \right) &\to 0 & \textrm{ as } |x|\to \infty \; . 
\end{aligned}
\end{equation}
\end{theorem}

\begin{remark}
  Hettlich proved these results for the derivatives of the far field operator
  in~\cite{hettlich1998frechet,hettlich1995frechet}.
  Those arguments can be adapted to obtain the derivatives of the
  values of the scattered field at a finite
  distance~\cite{borges2022multifrequency}.
\end{remark}

\begin{remark}
  Observe that the formula for $D_\Gamma \Fimp$ depends only on the normal
  component of the perturbation.
\end{remark}


For the $\lamch$ and $\lamabv$ models of the impedance function,
the impedance function depends on the curvature of the domain. The derivative of $\Fimp$ with
respect to $\Gamma$ can then be obtained by applying the chain rule

$$ \left [D_\Gamma \Fimp(\Gamma,\lambda(\Gamma)) \right ] \hb =
\left [\partial_\Gamma \Fimp(\Gamma,\lambda(\Gamma))\right ] \hb +
\left [\partial_\lambda \Fimp(\Gamma,\lambda(\Gamma)) \right ] \left [D_\Gamma \lambda(\Gamma)
  \right ] \hb \; ,$$
where the notations $\partial_\Gamma$ and $\partial_\lambda$ indicate the
Fr\'{e}chet derivative holding the other variable fixed.
The formula for $D_\Gamma \lambda(\Gamma)$ for both the $\lamch$
and $\lamabv$ models requires a formula for the Fr\'{e}chet derivative
of the curvature with respect to the boundary. 

Suppose that $\gammab:[0,L) \to \bbR^2$ is the parameterization
of a smooth curve, where $\gammab(s)=(x(s),y(s))$, and $H(\gammab)$ is
the signed curvature as defined in \cref{eq:signedcurve}.
Let $\hb$ be a sufficiently smooth $\bbR^2$-valued function on $[0,L)$.
Define a new curve
with the parameterization $\gammab_\eps(s) = \gammab(s) + \eps \hb(s)$.
The Fr\'{e}chet derivative $D_\Gamma H$ is then the linear operator
such that 

\begin{equation}
  \label{eq:Hder}
  \lim_{\epsilon \to 0} \frac{\left \| H(\gammab_\eps)
    - \left ( H(\gammab) + \left [ D_\Gamma H\right ] \hb \right ) \right \|}
      {\epsilon} = 0 
\end{equation}
for all such $\hb$.

The following direct formula for this derivative can be verified by
hand calculation:

\begin{proposition}
  Let $\hb$ be a normal perturbation of the curve, i.e.
  $\hb(t) = h(t)\nb(t)$ for a smooth function $h$. Then,
  \begin{equation}
    \label{eq:Hderdirect}
    \left [D_\Gamma H(\gammab) \right ] \hb = -H(\gammab)^2 h
    + \frac{x'x'' + y'y''}{(x'^2+y'^2)^2} h' -
    \frac{1}{x'^2+y'^2} h'' \; .
  \end{equation}
\end{proposition}

\section{Optimization methods}

\label{sec:optimization_methods}
This section describes details related to the optimization procedures
used to solve the inverse scattering problems. The continuation-in-frequency approach~\cite{chen1995recursive,chen1997inverse,bao2004inverse,bao2005inverse,bao2007inverse,bao2012shape,bao2015inverse,sini2012inverse,borges2015inverse,borges2017high,askham2023random} is the over-arching framework. This applies to
multi-frequency
data collected for a set of frequencies $\omega_1 < \omega_2 < \cdots < \omega_{N_k}$.
As explained in the introduction,
this approach begins with initial guesses, $\hat{\Gamma}_0$ and $\hat{\lambda}_0$, of
the domain boundary and impedance function, respectively, and does the iteration

\begin{equation}
  \label{eq:invprob_redux} 
  \left[\hat{\Gamma}_j,\hat{\lambda}_j \right]= \locargmin_{\Gamma \in S_\Gamma(\omega_j), \lambda \in S_\lambda(\omega_j)} \sum_{i=1}^{N_d} | \umeas_{\omega_j,\db_i} -
  \Fimp_{\omega_j,\db_i}(\Gamma,\lambda)|^2 \; ,
\end{equation}
where the notation $\locargmin$ indicates that the minimization problem is
solved using a standard local optimization
procedure based on gradient information and the initial guess
$\left [ \hat{\Gamma}_{j-1},\hat{\lambda}_{j-1} \right ]$.
The basic idea of the method is that the global minimizers at each frequency, subject
to appropriate choices of the constraint sets, should be sufficiently close
that a local solve works to find the global minimum at each step. See the
references above for more discussion.

\Cref{sec:opt_alt} provides an
overview of the alternating minimization approach we use to solve problems
of the form \cref{eq:invprob_redux}, which alternates between
minimization of the objective over the domain parameters and over the impedance
parameters.
We use two different methods for obtaining descent directions (in terms of
either the domain or impedance function parameters): Gauss-Newton
and steepest descent. Formulas for these in the inverse obstacle setting are
provided in \cref{sec:opt_descent_directions}. We discuss the strategies used
to minimize the objective function under the given constraints in the case of the
domain in \cref{sec:opt_filt} and in the case of the impedance function in
\cref{sec:projgrad}. 

\subsection{Alternating minimization}\label{sec:opt_alt}

We solve optimization problems of the form \cref{eq:invprob_redux}
by alternating minimization: we apply a step of an optimization algorithm
with respect to the domain boundary parameters and then a step of an optimization
algorithm with respect to the impedance function parameters.
This approach is classical~\cite{von1949rings,wright2015coordinate}
and has been previously applied in the
non-convex setting~\cite{bolte2014proximal,attouch2010proximal}. This is
particularly convenient because our methods for staying in the constraint
sets for the domain parameters and the impedance function parameters are
different. 

At frequency $\omega_j$, the initial guesses for the domain boundary and impedance
function are $\Gamma_{j}^{(0)} = \hat{\Gamma}_{j-1}$ and
$\lambda_{j}^{(0)} = \hat{\lambda}_{j-1}$, respectively. Then, we apply the following
alternating minimization iteration to find a local solution of \cref{eq:invprob_redux}:
\begin{itemize}
\item (Fix $\lambda$ and optimize over $\Gamma$)
  Let $\gammab:[0,L)\to \Gamma_{j}^{(\ell)}$ be a parameterization of
 $\Gamma_{j}^{(\ell-1)}$. Find the parameterization of a descent direction $\hb = h \nb$
 with respect to the domain boundary, where $h:[0,L)\to \bbR$ such that the curve $\Gamma_{j}^{(\ell)}$
 with the parameterization $\gammab + \hb:[0,L)\to \Gamma_{j}^{(\ell)}$ satisfies
$\Gamma_{j}^{(\ell)} \in S_\Gamma(\omega_j)$ and 
  \begin{equation}
    \label{eq:invprob_domain} 
    \sum_{i=1}^{N_d} | \umeas_{\omega_j,\db_i} -
    \Fimp_{\omega_j,\db_i}(\Gamma_{j}^{(\ell)},\lambda_{j}^{(\ell-1)})|^2 \leq  \sum_{i=1}^{N_d} | \umeas_{\omega_j,\db_i} -
    \Fimp_{\omega_j,\db_i}(\Gamma_{j}^{(\ell-1)},\lambda_{j}^{(\ell-1)})|^2 \; .  
  \end{equation}
  
\item (Fix $\Gamma$ and optimize over $\lambda$)
  Let $\gammab:[0,L)\to \Gamma_{j}^{(\ell)}$ be a parameterization of $\Gamma_{j}^{(\ell)}$.
    Find the parameterization of a descent direction $g:[0,L)\to \bbC$ with respect
      to the impedance
      function such that the impedance
      function $\lambda_{j}^{(\ell)}(s) = \lambda_{j}^{(\ell-1)}(s) + g(s)$ satisfies
    $\lambda_{j}^{(\ell)} \in S_\lambda(\omega_j)$ and   
    \begin{equation}
      \label{eq:invprob_imp} 
      \sum_{i=1}^{N_d} | \umeas_{\omega_j,\db_i} -
      \Fimp_{\omega_j,\db_i}(\Gamma_{j}^{(\ell)},\lambda_{j}^{(\ell)})|^2 \leq
        \sum_{i=1}^{N_d} | \umeas_{\omega_j,\db_i} -
    \Fimp_{\omega_j,\db_i}(\Gamma_{j}^{(\ell)},\lambda_{j}^{(\ell-1)})|^2    \; .
    \end{equation}
    
  \item Check stopping criteria.
\end{itemize}

Methods for finding appropriate descent directions and satisfying the constraints
are discussed in the following subsections. 

We impose a number of stopping criteria. We set a tolerance $\epsilon_R$ for the
relative residual

\begin{equation}
\label{eq:relres}
R_\omega(\hat{\Gamma},\hat{\lambda}) = \sqrt{\frac{\sum_{i=1}^{N_d} \left | \umeas_{\omega,\db_i} -
  \Fimp_{\omega,\db_i}(\hat{\Gamma},\hat{\lambda}) \right |^2}
{\sum_{i=1}^{N_d} \left | \umeas_{\omega,\db_i} \right |^2 } }\; ,
\end{equation}
such that we terminate if $R_\omega \leq \epsilon_R$. We try to detect stagnation
by monitoring for small relative changes in the norm of the domain parameters, impedance
parameters, and the relative residual. These tolerances are denoted by 
$\epsilon_{s,\Gamma}$, $\epsilon_{s,\lambda}$, and $\epsilon_{s,R}$, respectively.
We also impose a maximum number of iterations, $N_f$, of the alternating minimization
framework at each frequency.


\subsection{Obtaining descent directions}
\label{sec:opt_descent_directions}

The objective function in \cref{eq:invprob_redux} is in the form of a nonlinear
least squares problem. Here we briefly review two popular descent methods for such
objective functions: steepest descent and Gauss-Newton.

We review the real-parameter case; the complex-parameter case is similar.
Consider a multivariate function $F:\bbR^n \to \bbR$ defined by

\begin{equation*}
  F(\vb) = \frac{1}{2} \sum_{j=1}^m |f_j(\vb)-z_j|^2 = \frac{1}{2} \| \fb(\vb)-\zb\|^2 \; ,
\end{equation*}
where each $f_j:\bbR^n\to \bbC$ and $z_j \in \bbC$ for $j=1,\ldots,m$ and
$\fb(\vb)$ and $\zb$ are the vectors collecting these values. 
Recall that the Jacobian is the matrix
$\Jb(\vb)\in \bbC^{m\times n}$ with the entries $J_{ji}(\vb) = \partial_{x_i} f_j(\vb)$.

The gradient of $F$ is then

\begin{equation*}
  \nabla_v F = \begin{bmatrix} \real(\Jb(\vb)) \\ \imag(\Jb(\vb)) 
  \end{bmatrix}^T \begin{bmatrix} \real( \fb(\vb) - \zb ) \\
    \imag( \fb(\vb) - \zb ) 
  \end{bmatrix} \; .
\end{equation*}
Using the negative of the gradient as the search direction is known as
steepest descent minimization.

A reasonable step size, $d$, in the direction of the negative of the gradient can be
obtained by minimizing the approximation

\begin{align*}
  F(\vb -d\nabla_vF) &\approx \frac{1}{2} \| \fb(\vb) - d \Jb(\vb) \nabla_v F - \zb\|^2 \\
  &= \frac{1}{2} \|\fb(\vb)-\zb\|^2 - d \| \nabla_v F \|^2  + \frac{d^2}{2} \|J(\vb) \nabla_vF\|^2 \; ,
\end{align*}
as a function of $d$. The minimum occurs for

\begin{equation}
  \label{eq:cauchypoint}
  d = \frac{ \|\nabla_v F\|^2 }{ \|J(\vb) \nabla_v F\|^2 } \; .
\end{equation}
This is the so-called Cauchy point for the linearization. We use it as an initial
guess for the step size when the direction is based on the gradient.

The Gauss-Newton method is based on the linearization $\fb(\vb + \Delta \vb)
\approx \fb(\vb) + \Jb(\vb) \Delta \vb$. The minimizer $\Delta \vb$ of
the approximation

\begin{equation*}
F(\vb + \Delta \vb) \approx \| \fb(\vb) - \zb + \Jb(\vb) \Delta \vb \|^2 
\end{equation*}
is then the linear least squares solution of the system
\begin{equation*}
\begin{bmatrix} \real(\Jb(\vb)) \\ \imag(\Jb(\vb)) 
  \end{bmatrix} \Delta \vb \approx -\begin{bmatrix} \real( \fb(\vb) - \zb ) \\
    \imag( \fb(\vb) - \zb ) 
  \end{bmatrix}   \; ,
\end{equation*}
which can be computed by direct methods. The Gauss-Newton step has a built-in
notion of the appropriate step size, so the initial guess for $d$ is 1 for this
descent direction. 

\paragraph{Specifics for boundary curve optimization}
Suppose that at frequency $\omega$ we have the approximate
boundary curve, $\Gamma$, of length $L$ and the approximate impedance
function $\lambda$. Suppose that this curve has the arc-length
parameterization
$\gammab:[0,L) \to \bbR^2$. As noted above, we consider only normal perturbations of this curve
of the form $\hb = h \nb$, where $\nb(s)$ is the normal to the curve at $\gammab(s)$.
Let $N_\gamma(\omega)$ be given and $h_\wb$ be parameterized as in \cref{eq:hupdate}, i.e.
as a real Fourier series, with
parameters $\wb = (a_{h,0},a_{h,1},\ldots,a_{h,N_\gamma},b_{h,0},b_{h,1},\ldots,b_{h,N_\gamma}) \in \bbR^{2N_\gamma+1}$.
Then, we consider minimizing the objective function 

\begin{equation*}
F_\Gamma(\wb) = \sum_{i=1}^{N_d} | \umeas_{\omega,\db_i} -
\Fimp_{\omega,\db_i}(\Gamma(\wb),\lambda)|^2 \; ,
\end{equation*}
where the curve $\Gamma(\wb)$ is the curve with the (not necessarily arc-length)
parameterization $\gammab_\wb:[0,L)\to \bbR^2$ defined by 
$\gammab_\wb(s) = \gammab(s) + h_\wb(s) \nb(s)$.

The values of

$$ \partial_{\text w_j} \Fimp_{\omega,\db_i}(\Gamma(\wb),\lambda) \in \bbC^{N_r}$$
can be obtained using the Fr\'{e}chet derivative formulae of \cref{sec:frechetderivative}
and applying the chain rule when appropriate. In particular,
\begingroup
\setlength{\tabcolsep}{10pt} 
\renewcommand{\arraystretch}{2} 
\begin{equation*}
  \partial_{\text w_j} \Fimp_{\omega,\db_i}(\Gamma(\wb),\lambda)
  = \left \{ \begin{array}{rcl}
    \left [D_\Gamma \Fimp_{\omega,\db_i}(\Gamma,\lambda) \right ] \phi_j \nb
    & \textrm{ if } & j \leq N_\gamma + 1 \\     
    \left [D_\Gamma \Fimp_{\omega,\db_i}(\Gamma,\lambda) \right ] \psi_j \nb
    & \textrm{ if } & j > N_\gamma + 1 
  \end{array} \right. \; ,
\end{equation*}
\endgroup
where $\phi_j(s) = \cos(2\pi(j-1)s/L)$ and $\psi_j(s) = \sin(2\pi(j-N_\gamma-1)s/L)$.
Likewise, the case in which
$\lambda$ depends on $\Gamma$ can be treated; we omit the details for
the sake of brevity.

\begin{remark}
Observe that the Jacobian matrix has $\Jb(\wb) \in \bbC^{N_r N_i \times (2N_\gamma + 1)}$.
We must solve a PDE for each incident direction $\db_i$ and each entry of $\wb$
to fill the matrix. See \cref{sec:pdenumerics} for a brief discussion of
computational costs.
\end{remark}

\paragraph{Specifics for impedance function optimization}

Suppose that at frequency $\omega$ we have the approximate
boundary curve, $\Gamma$, of length $L$. We consider three
objective functions, depending on the impedance function model:

\begin{align*}
F_{\lamfs}(\cb) &= \sum_{i=1}^{N_d} | \umeas_{\omega,\db_i} -
\Fimp_{\omega,\db_i}(\Gamma,\lamfs[\cb])|^2 \; , \\
F_{\lamch}(\alphab) &= \sum_{i=1}^{N_d} | \umeas_{\omega,\db_i} -
\Fimp_{\omega,\db_i}(\Gamma,\lamch[\alphab])|^2 \; , \\
F_{\lamabv}(\betab) &= \sum_{i=1}^{N_d} | \umeas_{\omega,\db_i} -
\Fimp_{\omega,\db_i}(\Gamma,\lamabv[\betab])|^2 \; . \\
\end{align*}

As above, the entries of the Jacobian can be obtained using
the Fr\'{e}chet derivative formulae of \cref{sec:frechetderivative}
and applying the chain rule when appropriate. For example,
\begingroup
\setlength{\tabcolsep}{10pt} 
\renewcommand{\arraystretch}{2} 
\begin{equation*}
  \partial_{\alpha_j} \Fimp_{\omega,\db_i}(\Gamma,\lamch[\alphab])
  = \left \{ \begin{array}{rcl}
    \left [D_\lambda \Fimp_{\omega,\db_i}(\Gamma,\lamch[\alphab]) \right ] 1 
    & \textrm{ if } & j = 1 \\
    \left [D_\lambda \Fimp_{\omega,\db_i}(\Gamma,\lamch[\alphab]) \right ] H 
    & \textrm{ if } & j = 2 
  \end{array} \right. \; ,
\end{equation*}
\endgroup
where $1$ denotes the function that is constant 1 and $H$ is the signed
curvature function on $\Gamma$. The $\lamfs$ and $\lamabv$
cases can be treated in similar fashion; we omit the details for
the sake of brevity.

\subsection{Filtering methods for the domain}
\label{sec:opt_filt}

Following the notation of the previous section, let $\wb$ be the
Fourier series coefficients defining an update $h_{\wb}$. We find coefficients
$\wb_\gn$ and $\wb_\sd$ corresponding to Gauss-Newton and steepest descent directions,
respectively, and set $d$ to be the step length for steepest descent determined using the
Cauchy point formula \cref{eq:cauchypoint}. We set $h^\gn = h_{\wb_\gn}$ and
$h^\sd = h_{d\wb_\sd}$. For either update, it is not guaranteed that the new curve,
$\Gamma(\wb_\gn)$ with parameterization $\gammab_\gn = \gammab + h^\gn \nb$ or
$\Gamma(d\wb_\sd)$ with parameterization $\gammab_\sd = \gammab + h^\sd \nb$,
is in the constraint set, $S_\Gamma(\omega)$. We follow two strategies for modifying these
updates~\cite{askham2023random}: step-length filtering and Gaussian filtering.

Step-length filtering is the basic strategy of shortening the step in
the given direction until the constraints are met and the residual is
non-increasing. Let $\eta_\filt > 1$ and $N_{\filt}$, a positive integer,
be the filtering parameters. For any update $h$ with parameters $\wb$,
the filtered update is then $h/\eta_\filt^{n_\filt}$, where $n_\filt$ is the smallest
$\ell$ with $0\leq \ell \leq N_\filt$ such that 
\begin{equation*}
  \Gamma(\wb/\eta_\filt^\ell) \in S_\Gamma(\omega)
  \quad \textrm{and} \quad F_\Gamma(\wb/\eta_\filt^\ell) \leq F_\Gamma(\zerob) \; ,
\end{equation*}
with the convention that the filtered step is of size 0 if there is no $\ell$ satisfying
the conditions.
We denote the step-length filtered steps for each direction by $h^{\gn,\slf}$ and
$h^{\sd,\slf}$.

In practice, step-length filtering results in overly short steps when the update
induces self-intersections or high curvature. The idea of Gaussian filtering is to
attenuate high frequency components of the update in order to more efficiently avoid
these geometric issues~\cite{askham2023random}.

Let $\wb$ be the parameters of a given update. Recall that the first $N_\gamma+1$ entries
correspond to cosines of increasing frequency and the last $N_\gamma$ entries to sines
of increasing frequency. Let $\sigma_\filt > 0$ be a given filtering parameter determining
the width of a Gaussian and define the diagonal linear operator $G_{\sigma}$ by 
\begingroup
\setlength{\tabcolsep}{10pt} 
\renewcommand{\arraystretch}{2} 
\begin{equation*}
G_{\sigma}(\eb_m)
= \left \{ \begin{array}{lcl}
			\exp\left(-\frac{(m-1)^2}{\sigma^2N_h^2}\right) \eb_m & \textrm{ if } & 1\leq m\leq N_h+1,\\
            \exp\left(-\frac{(m-N_h-1)^2}{\sigma^2N_h^2}\right) \eb_m & \text{ if } & N_h+2\leq m\leq 2N_h+1. 
            \end{array} \right. \; 
\end{equation*}
\endgroup
The filtered update is then $h_{\wb_\gf}$, where $\wb_\gf = G_{\sigma_\filt^{n_\filt}}\wb$ and $n_\filt$
is the smallest $\ell$ with $0\leq \ell \leq N_\filt$ such that 
\begin{equation*}
  \Gamma(G_{\sigma_\filt^{\ell}}\wb) \in S_\Gamma(\omega) \quad \textrm{and} \quad
  F_\Gamma(G_{\sigma_\filt^{\ell}}\wb) \leq F_\Gamma(\zerob) \; ,
\end{equation*}
with the convention that the filtered step is of size 0 if there is no $\ell$ satisfying the
conditions. 
We denote the Gaussian filtered steps for each direction by $h^{\gn,\gf}$ and
$h^{\sd,\gf}$. 

We have found that the best performance is obtained by attempting all strategies.
At each step, we select $h$ to be the filtered step out of the options
$h^{\gn,\slf}$, $h^{\gn,\gf}$, $h^{\sd,\slf}$, and $h^{\sd,\gf}$ for which the updated
curve results in the minimum residual.

The updated curve is defined by the parameterization $\gammab + h \nb$, which is
not necessarily an arc-length parameterization.
We re-parameterize the domain in arc-length using the algorithm described
in~\cite{beylkin2014fitting}, which is based on fitting a bandlimited curve
to equispaced (in parameter space, not necessarily arc-length)
samples of the new parameterization.

\subsection{Projection methods for the impedance function}
\label{sec:projgrad}

The constraint sets for the $\lamch$ and $\lamabv$ models are
closed and convex, so the corresponding constrained minimization
problem can be treated by a projected gradient method \cite{zarantonello1971projections,conn2000trust}. 
In particular, the impedance optimization problems are of the form

$$ \min_{\vb \in \Ccal} F(\vb) \; ,$$
where $\mathcal{C}$ is a closed, convex set. The simple iteration

\begin{equation}
  \label{eq:projiter}
  \vb_{j+1} = \proj_{\Ccal} (\vb_j - d_j \nabla_v F(\vb_j)) \; ,
\end{equation}
where $d_j$ is a step-length parameter and 

$$ \proj_{\Ccal} (\yb) = \argmin_{\xb \in \Ccal} \|\yb - \xb\| \; ,$$
is known as a projected gradient method. The iteration is known to converge
under mild assumptions on $f$~\cite{gafni1982convergence, conn2000trust}.

The convex constraint set for $\lamch$ is given by

$$ \Ccalch = \{ \vb \in \bbC^{2} : \imag(\vb) \leq \zerob \} \; ,$$
where the imaginary part and inequality are interpreted component-wise.
The projection operator for this set is simple to compute:

$$ \proj_{\Ccalch} (\yb) = \real(\yb) + \imath \min(\imag(\yb),\zerob) \; ,$$
where, likewise, the minimum and real part operators are interpreted
component-wise. The convex constraint set for $\lamabv$ is given by

$$ \Ccalabv = \{ \vb \in \bbR^{2} : \vb \geq \zerob \} \; .$$
The projection operator for this set is also simple to compute:

$$ \proj_{\Ccalabv} (\yb) = \max(\yb,\zerob) \; .$$

We use a similar strategy to the step length filtering strategy for
the domain to select the step length. For example, in the $\lamabv$
model the updated parameters would be

$$ \proj_{\Ccalabv} \left (\betab - d \nabla_\beta F_\lamabv(\betab)/\eta_\filt^{n_\filt} \right ) \; ,$$
where $d$ is the original step length determined by the Cauchy point
formula \cref{eq:cauchypoint} and $n_\filt$ is the smallest $\ell$ with
$0\leq \ell \leq N_\filt$ such that 

\begin{equation*}
  F_\lamabv \left ( \proj_{\Ccalabv} \left (\betab - d \nabla_\beta F_\lamabv(\betab)/\eta_\filt^{\ell} \right ) \right )
  \leq F_\lamabv(\betab) \; ,
\end{equation*}
with the convention that the step is size 0 if there is no
$\ell$ satisfying the conditions.

\begin{remark}
  While we did not use a Gauss-Newton descent direction for the impedance
  function optimization in the examples below, it can be done and
  there is one interesting phenomenon to note. If the curvature of the
  domain is nearly constant, say for a circular domain, the Jacobian
  matrix for the curvature dependent models is highly ill-conditioned.
  We found that one remedy was to compute the condition number of the Jacobian
  and revert to a constant impedance model if the condition number was above
  some threshold.
\end{remark}

\section{Numerical results}\label{sec:numerical}

The numerical experiments in this section explore the
suitability of impedance models for dissipative transmission
problems. We describe common test settings here. Tests with
data generated by an impedance model are in \cref{sec:impdata}
and tests with data generated by a transmission model are
in \cref{sec:transdata}.

\paragraph{Reproducibility}
The scripts used to generate these results~\cite{travis_askham_2023_10214846}
and the data associated with the figures~\cite{askham_2023_10214866}
are publicly available and archived. 

\paragraph{Test data parameters}
Each test has a range of pulsation frequencies, which are
determined by the value $\ktwomax$, taken to be an
integer. In particular, the values of $\omega$ are of the form 

$$ \omega_j = c_2 \left( 1 + (j-1)/2 \right ) $$
where $j = 1,\ldots,2\ktwomax+1$. Thus, there are $N_k=2\ktwomax+1$
frequencies and the exterior wave speeds range from $1$ to $\ktwomax$,
with a spacing of $1/2$. The value of $\ktwomax$ is selected so
that the obstacle can be reasonably well reconstructed from data
at the highest frequency.

For most of the examples, we assume best-case data in the sense
that there are sufficiently many incident directions and receptor
locations to obtain a good reconstruction and that the data
are full aperture. In particular, unless otherwise noted, the
data are collected at all receptors for each incident angle, 
the number of incident directions and receptors is set to
$N_d(\omega) = N_r(\omega) = \lfloor 10 \omega/c_2 \rfloor$,
the incident directions are set to $\db_i = (\cos(\theta_i),\sin(\theta_i))$
for $i = 1,\ldots,N_d(\omega)$
with $\theta_i = 2\pi i/N_d(\omega)$, and the receptors are located
at the points $\rb_i = (r\cos(\theta_i),r\sin(\theta_i))$
for $i=1,\ldots,N_r(\omega)$ with $r=10$.

One of our primary interests below is the performance of the
model as the dissipation, $\delta$, varies. As a reference
dissipation value, we use $\delta_0 = \sqrt{3} \ktwomax$, which
is suggested in~\cite{antoine2005construction} as a
value of dissipation where the asymptotic model should
become accurate for the highest pulsation in the data.
The remaining physical parameters are fixed
as $c_1=0.5$, $c_2=1.0$, $\rho_1=1.2$, and $\rho_2=0.7$,
so that $c_r = 0.5$ and $\rho_r \approx 1.7$. This corresponds
to the obstacle having a denser material and a higher acoustic
wave speed. We found that
the results were similar if instead $c_r < 1$ and $\rho_r > 1$,
though we did not explore the extremes.

\paragraph{Measures of error} 
When applying continuation-in-frequency, we obtain a sequence of
approximations of the domain, $\{\hat{\Gamma}_j\}_{j=1}^{N_k}$, and
impedance functions, $\{\hat{\lambda}_j\}_{j=1}^{N_k}$, corresponding
to the solution obtained at the frequencies $\{\omega_j\}_{j=1}^{N_k}$.
We consider two quantitative error measurements for a given reconstruction
of the domain and impedance function $\hat{\Gamma},\hat{\lambda}$. 
The first is the relative residual at a given frequency; see~\cref{eq:relres}.
Below, we always plot $R_{\omega_j}(\hat{\Gamma}_j,\hat{\lambda}_j)$.
The second error measurement is the relative area of the
symmetric difference between the interior of the
true obstacle, $\Omega_\star$, and the interior of the
recovered obstacle, $\hat{\Omega}$ (which has boundary
$\hat{\Gamma}$):

\begin{equation}
  \label{eq:area}
  E(\hat{\Gamma}) = \frac{\area(\Omega_\star \setminus \hat{\Omega})
    + \area(\hat{\Omega} \setminus \Omega_\star)}{\area(\Omega_\star)} \; .
\end{equation}
We approximate the areas in \cref{eq:area} based on the polygons
defined by the boundary nodes of $\Omega_\star$ and $\hat{\Omega}$.
This is an inherently low order approximation but appears to
provide at least 2 digits of precision in the examples below.

\paragraph{Discretization, regularization, and optimization method
  parameters}

The PDEs are discretized using high order boundary integral equation
methods as described in \cref{sec:bieform}.
We discretize the boundary curve using approximately 10 points per wavelength
by setting $N(\omega) = \max( \lceil 5 \omega L/(c_2\pi) \rceil, 300 )$,
where $L$ is the length of the current approximation of the domain.
This ensures that the forward map, $\Fimp$, and its Fr\'{e}chet derivatives
are evaluated accurately.

In most of the examples, the data are generated using the forward model with
transmission boundary conditions and the inverse problem is solved using the model 
with impedance boundary conditions, so that ``inverse crimes'' are avoided.
However, this is not the case in
\cref{sec:impdata}. An inverse crime is still avoided there in that
the data are generated using approximately 20 points per wavelength whereas
the inverse problem uses 10,
the non-airplane domains are not originally parameterized in arc-length whereas
the inverse solver constantly resamples in arc-length, and 
the true boundary curve is generally not contained in the constraint
set at the lowest frequencies.

The search direction for the domain boundary update, $h$, is regularized
by limiting the length of its real Fourier series. This length is $2N_\gamma(\omega)+1$
and $N_\gamma(\omega)$ should be proportional to $\omega$ so that $h$ does
not contain features smaller than half of a wavelength. We set
$N_\gamma(\omega) = \lfloor \omega L/(c_2 \pi) \rfloor$, where $L$ is
the length of the current approximation of the domain.
The constraint set $S_\Gamma(\omega)$ is determined by the amount of
curvature regularization used. The idea is that the amount of elastic
energy contained in the first $M(\omega)$ modes of the curvature
should exceed some proportion of the total elastic energy, i.e.
$\mathcal{E}^{M(\omega)} \geq \ccurv \mathcal{E}(\Gamma)$; cf. \cref{eq:em}.
We set $\ccurv = 0.9$ and $M(\omega) = N_\gamma(\omega)$.
The $\lamfs$ model is also regularized by limiting the length of its
Fourier series. This length is $2N_c(\omega)+1$ and we set
$N_c(\omega) = N_\gamma(\omega)/2$ following the advice of
\cite{borges2022multifrequency}.

For detailed descriptions of the optimization parameters, refer to
\cref{sec:optimization_methods}. Alternating minimization is performed
with a maximum of $N_f=40$ iterations, with tolerance $\epsilon_R = 10^{-4}$
and stagnation parameters $\epsilon_{s,\Gamma} = \epsilon_{s,\lambda} = \epsilon_{s,R} = 10^{-4}$.
We set a maximum of $N_\filt = 3$ filtrations per step. Step length filtering is performed
with $\eta_\filt = 8$ and Gaussian filtering is performed with
$\sigma_\filt = 10^{-1}$.

\subsection{Experiments with impedance data}

\label{sec:impdata}

\begin{figure}
  \centering
  \includegraphics[width=1\textwidth]{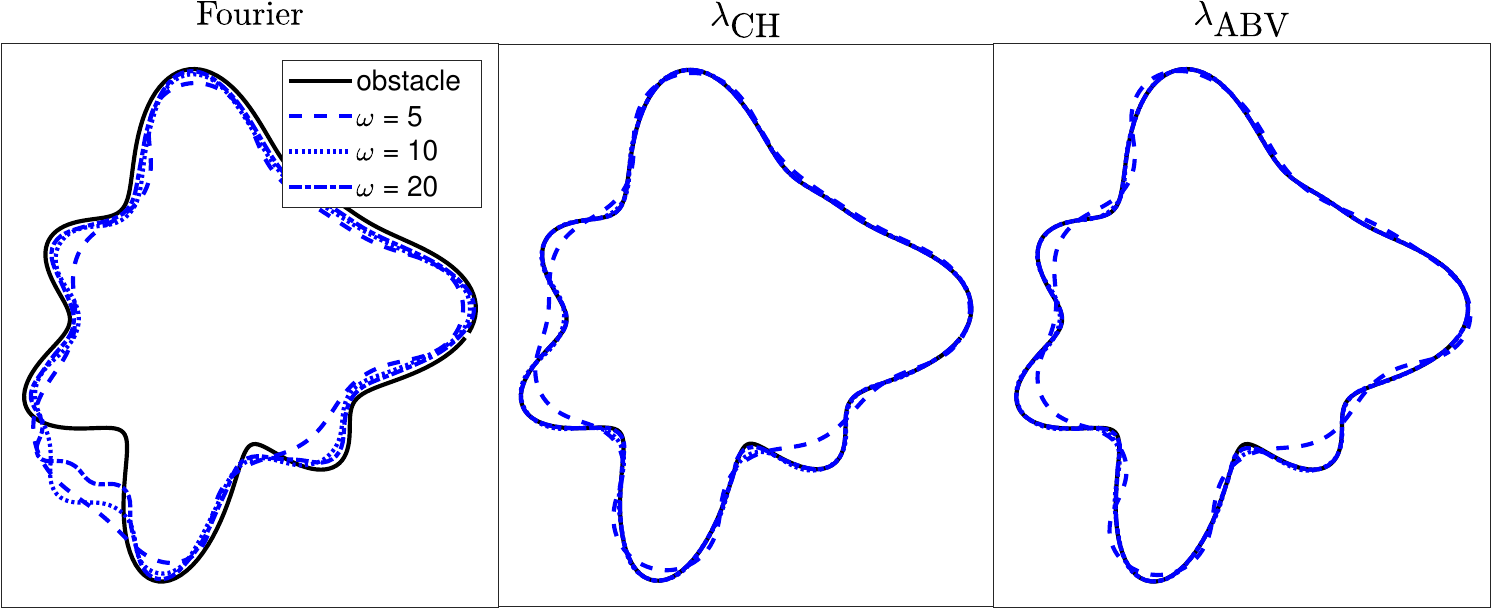}
    
  \vspace{.1cm}
  
  \includegraphics[width=1\textwidth]{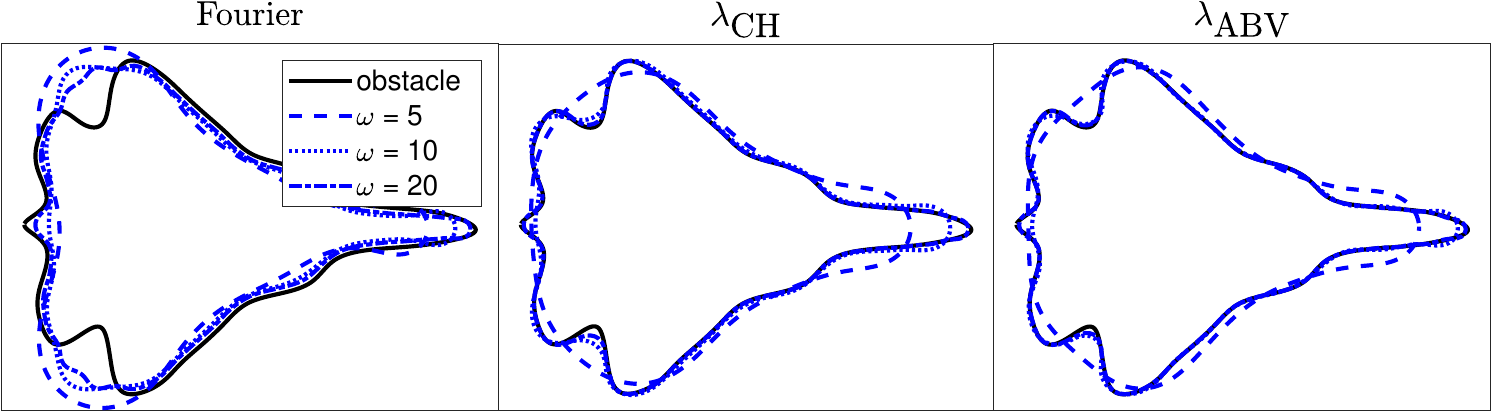}
  
  \vspace{.1cm}
    
  \includegraphics[width=1\textwidth]{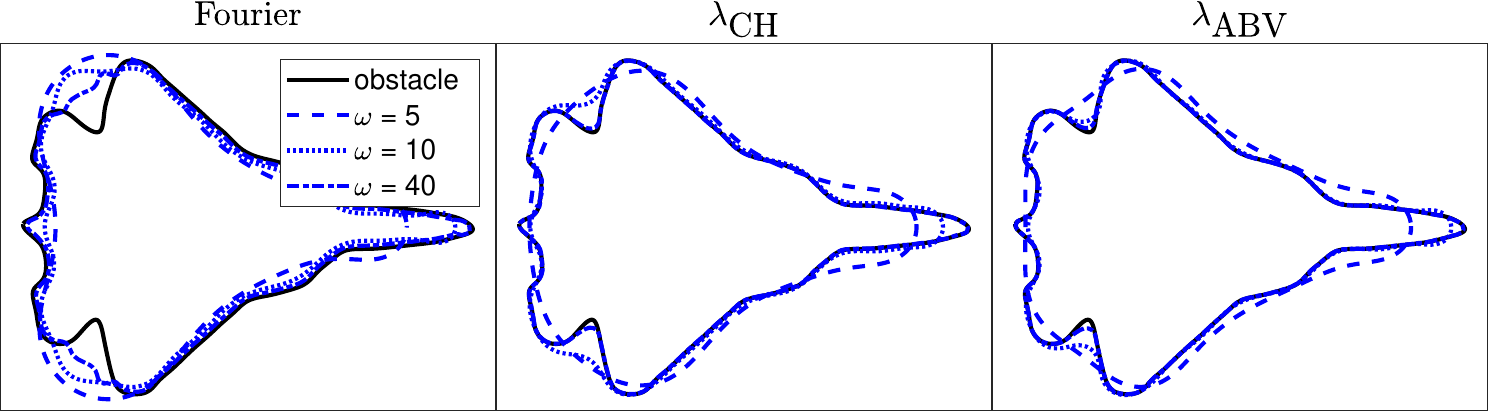}
    
  \vspace{.1cm}
  
  \includegraphics[width=1\textwidth]{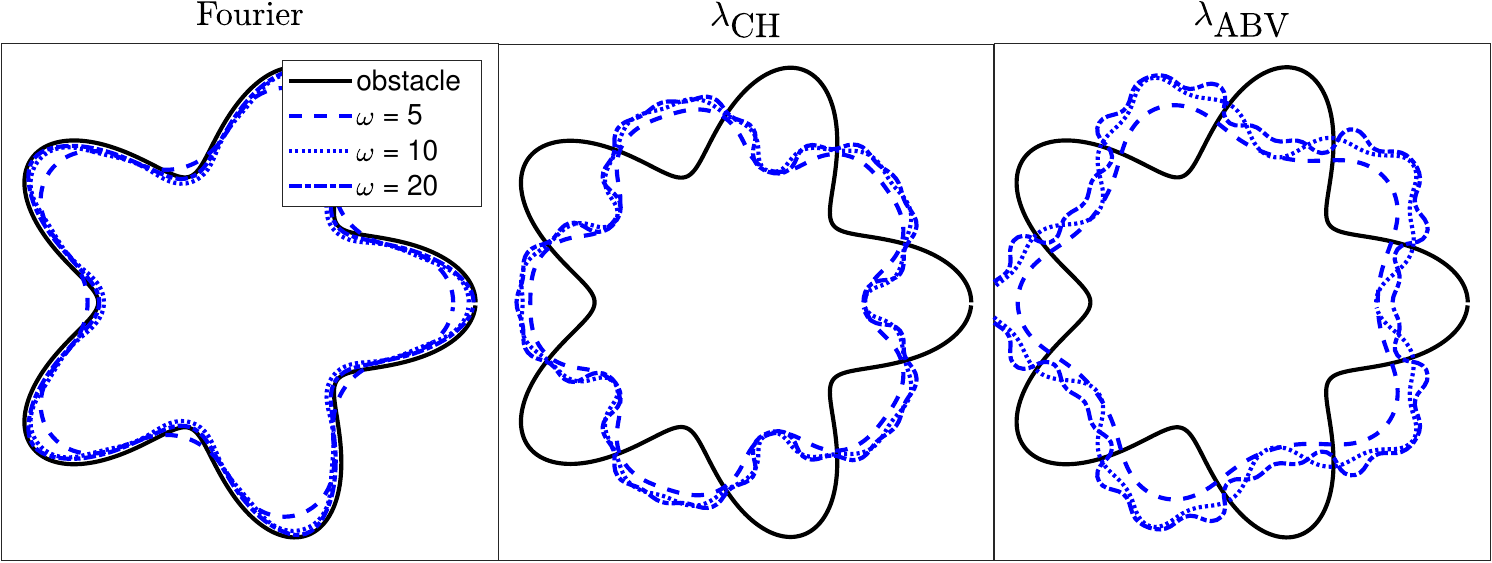}
    
  \vspace{.1cm}
  
  \caption{Reconstructions of the domains obtained for the experiments
    with impedance data of \cref{sec:impdata}. The top row corresponds to
    a ``random'' obstacle, the second row a smooth ``airplane''
    obstacle,
    the third row a more detailed ``airplane'' obstacle,
    and the bottom row a
    ``starfish'' obstacle. Each column is labeled by the impedance
    model used in the inverse problem.}
  \label{fig:diff_imp_models}
\end{figure}

Here we compare the performance of the three impedance models,
$\lamfs$, $\lamch$, and $\lamabv$, for impedance data generated
by the model of~\cite{antoine2001high}, i.e. generated as the
solution of \cref{eq:imppde} with the impedance
model \cref{eq:abvorder1} and the dissipation set as
$\delta = \delta_0 = \sqrt{3}\ktwomax$.
We consider 4 different
obstacle shapes. In \cref{fig:diff_imp_models}, the first
corresponds to a smooth obstacle shape parameterized as
$(r(t)\cos(t),r(t)\sin(t))$ with $r(t)$ given as constant perturbed
by a random Fourier series ($\ktwomax=20$), the second row is a smooth airplane-like
obstacle ($\ktwomax=30$), the third row is a more detailed airplane-like
obstacle  ($\ktwomax=40$), and
the fourth row is a smooth starfish-like obstacle ($\ktwomax=30$). 

The $\lamfs$ model performs poorly, even on the smooth, random
shape. In general, the $\lamfs$ model appears to get caught in
local minima. This is in contrast to previous studies, e.g.
\cite{borges2022multifrequency}. We believe that the important distinction
of the dissipative model setting is that we must allow for $\lambda$
to be a general, complex-valued function, whereas prior studies could
constrain $\lamfs$ to be real.

The constrained, curvature-dependent models, $\lamch$ and $\lamabv$, perform
much
better on the first three obstacles. To the eye, these are capable
of recovering sharp features and non-convex features
of the more detailed airplane. The $\lamabv$ model appears to
have a slight edge over the $\lamch$ model, which can be seen most
clearly by comparing the results obtained for the two at pulsation
$\omega = 10$.

\begin{figure}
  \centering
  \includegraphics[width=0.6\textwidth]{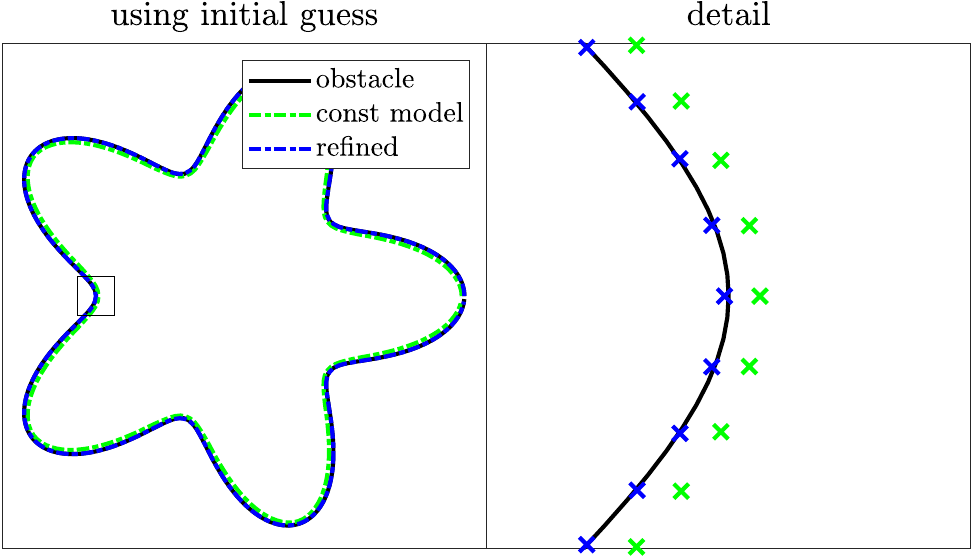}
  \caption{Using a different initializer for the starfish domain in
    the impedance data experiment of \cref{sec:impdata}. Here we
    show an initial guess obtained from solving the inverse problem with
    a constant impedance model and the refined inverse problem solution
    obtained by the $\lamabv$ model with this initial guess.}
  \label{fig:const_model_initial_guess}
\end{figure}

The constrained, curvature-dependent models fail on the
starfish-like obstacle. These models get caught in a local
minimum that has the same symmetry as the domain.
Based on the results of the random smooth domain, which
is similarly constructed, and the airplane domains, which have
reflective symmetry, it appears that this phenomenon is specific
to this type of starfish-like domain in which the distance from
the origin is correlated with the curvature of the domain.
The $\lamfs$ model avoids this issue because we select $N_c(\omega)$ to
be $N_\gamma(\omega)/2$, which forces the optimization problem to
resolve the domain first. In contrast, the curvature-dependent
models can have Fourier content that is similar to the obstacle.

While this failure mode is worthy of note, we believe it to be
unlikely to appear in applications. It can also be mitigated by
first solving the inverse problem with a constant impedance model
and then using the
resulting approximation of the obstacle as an initial guess to then
find a curvature-dependent model; see \cref{fig:const_model_initial_guess}
for an example of this,
where the approximation of the obstacle obtained for the constant
impedance model is quite good and the curvature-dependent
$\lamabv$ model succeeds with this initial guess. Because of the
superior performance of $\lamabv$ on most
of these examples, we use that model for the remaining examples,
unless otherwise noted. 

\subsection{Experiments with transmission data}

\label{sec:transdata}

The results of this section are for data generated by the transmission
problem, \cref{eq:transpde},
for the more detailed airplane-like obstacle. In all examples in this
section, $\ktwomax=40$.
The model used for the inverse problems is always
the impedance model, \cref{eq:imppde}, with the
$\lamabv$ model for the impedance function, unless otherwise noted. 

\subsubsection{Effect of lowering dissipation}

\label{sec:lowerdelta}

\begin{figure}
  \centering
  \includegraphics[width=1\textwidth]{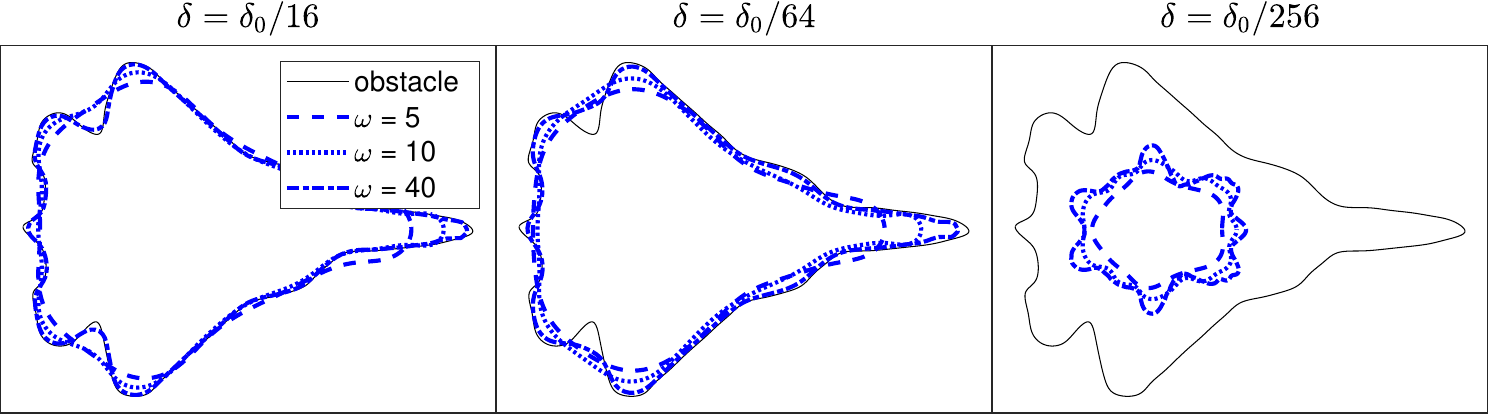}
  
  \vspace{.5cm}
  
  \includegraphics[width=0.8\textwidth]{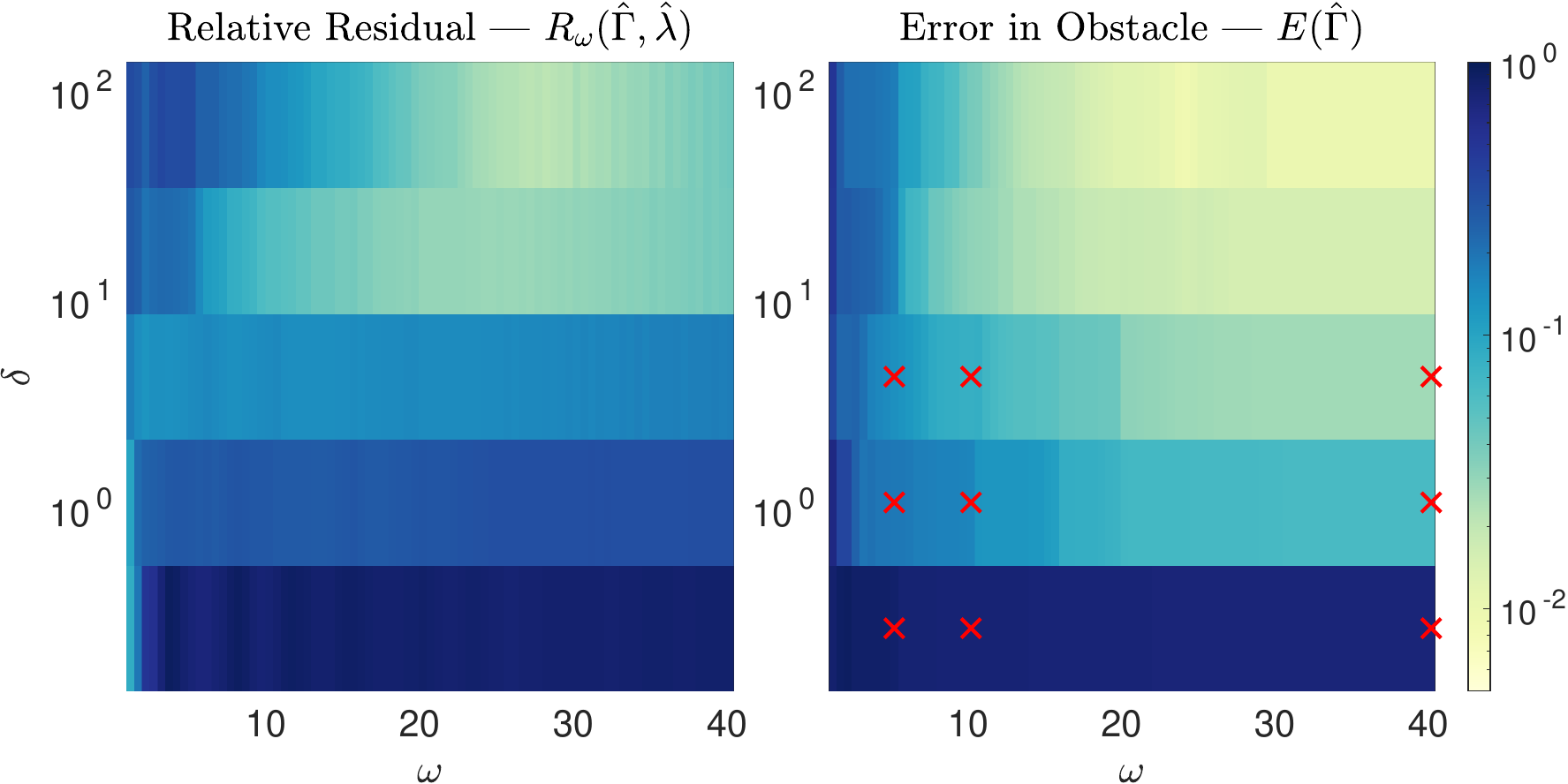}
  \caption{Experiment of \cref{sec:lowerdelta} with transmission data
    using the $\lamabv$ impedance model for the inverse problem.
    The top row shows reconstructions obtained for different
    values of the pulsation, $\omega$, and dissipation, $\delta$. The bottom row has plots
    of the error measures as a function of $\omega$ and $\delta$.
    The values of $(\omega,\delta)$ for the reconstructions in the top
    row are marked in red in the error plot.}
    \label{fig:resanderr_varydeltafreq}
\end{figure}

The $\lamabv$ model is known to be more accurate for higher
levels of dissipation, with $\delta_0 = \sqrt{3} \ktwomax$
providing good quantitative agreement in the scattered data
of the impedance and transmission models~\cite{antoine2005construction}.

In this experiment, we generate transmission scattering data
for the more detailed plane with $\delta = \delta_0/4^j$ for
$j=0,1,2,3,4$. We then solve the inverse problem with the
$\lamabv$ model for each data set. The results are
shown in \cref{fig:resanderr_varydeltafreq}. 

We find that a reasonably close reconstruction of the obstacle
is obtained for dissipation as low as $\delta=\delta_0/64$.
For $\delta=\delta_0/256$, the reconstruction fails completely.

\begin{figure}
  \centering
  \includegraphics[width=1\textwidth]{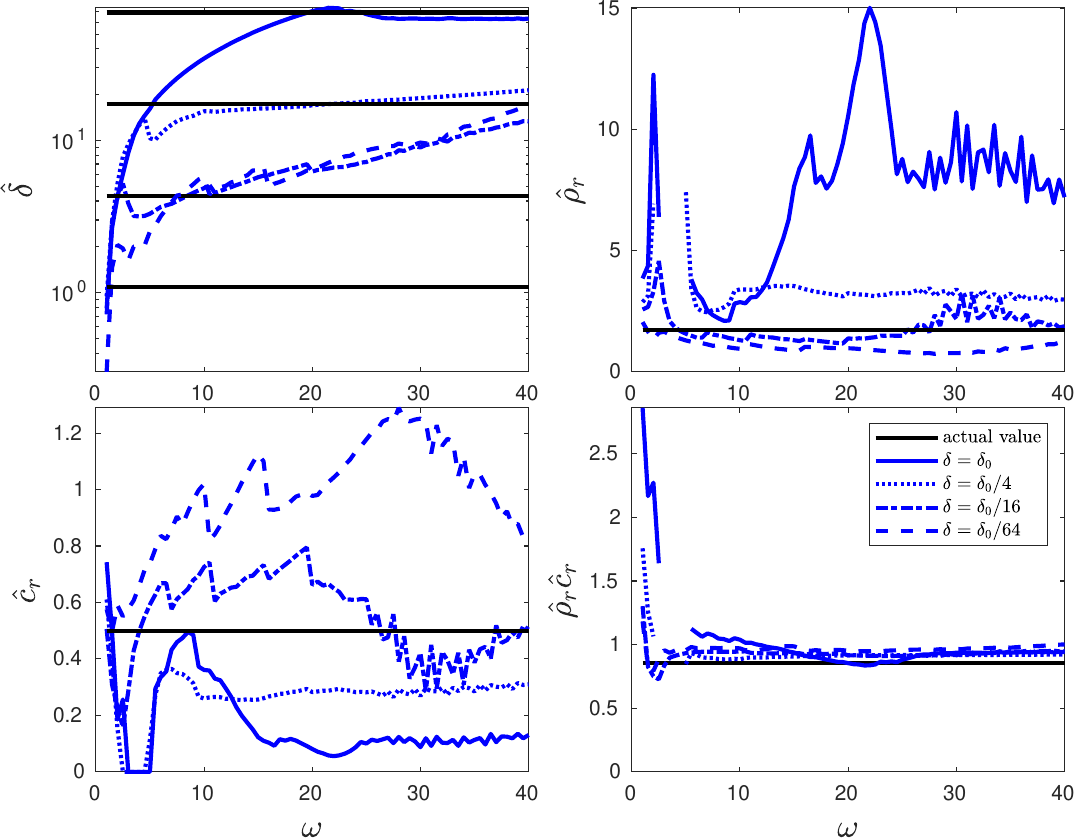}
  \caption{Experiment of \cref{sec:lowerdelta} with transmission data
    using the $\lamabv$ impedance model for the inverse problem. Each plot
    shows the recovery of a physical parameter for various values of
    the dissipation, $\delta$, as a function of the pulsation, $\omega$.}
    \label{fig:recovered_parameters_varydeltafreq}
\end{figure}

Given the parameters $\hat{\betab}$ recovered by solving the
inverse problem with the $\lamabv$ model, we can recover
approximations of the physical parameters, denoted
$\hat{\rho}_r$, $\hat{c}_r$, $\hat{\delta}$, from the
relations \cref{eq:betatophys}.
\Cref{fig:recovered_parameters_varydeltafreq} plots 
$\hat{\delta}$, $\hat{\rho}_r$, $\hat{c}_r$ and
$\hat{\rho}_r\hat{c}_r$ as a function of the
pulsation for various values of the dissipation, $\delta$.
Because the obstacle recovery improves as $\omega$
increases, it is expected that the recovered parameters should
be better for higher values of the pulsation.
The value of $\hat{\delta}$ is reasonably 
accurate for higher levels of dissipation, but it
gives an overestimate for lower levels of dissipation.
In most cases, the recovered values $\hat{\rho}_r$ and
$\hat{c}_r$ are not particularly accurate. However, the product
$\hat{\rho}_r\hat{c}_r$ is relatively accurate in most regimes. We suspect
that this is because the product $\hat{\rho}_r\hat{c}_r$ appears in the
constant term of the $\lamabv$ model, whereas only $\hat{\rho}_r$ appears
in the curvature term; cf. \cref{eq:abvorder1}. 

\subsubsection{Experiment with noise}

\label{sec:withnoise}

\begin{figure}
  \centering
  \includegraphics[width=1\textwidth]{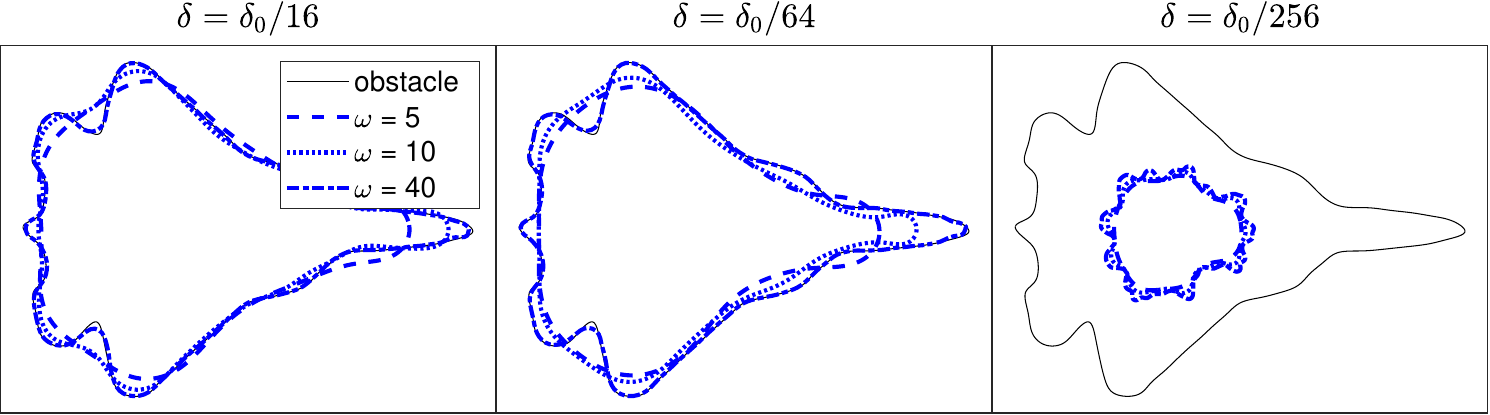}
  
  \vspace{.5cm}
    
  \includegraphics[width=0.8\textwidth]{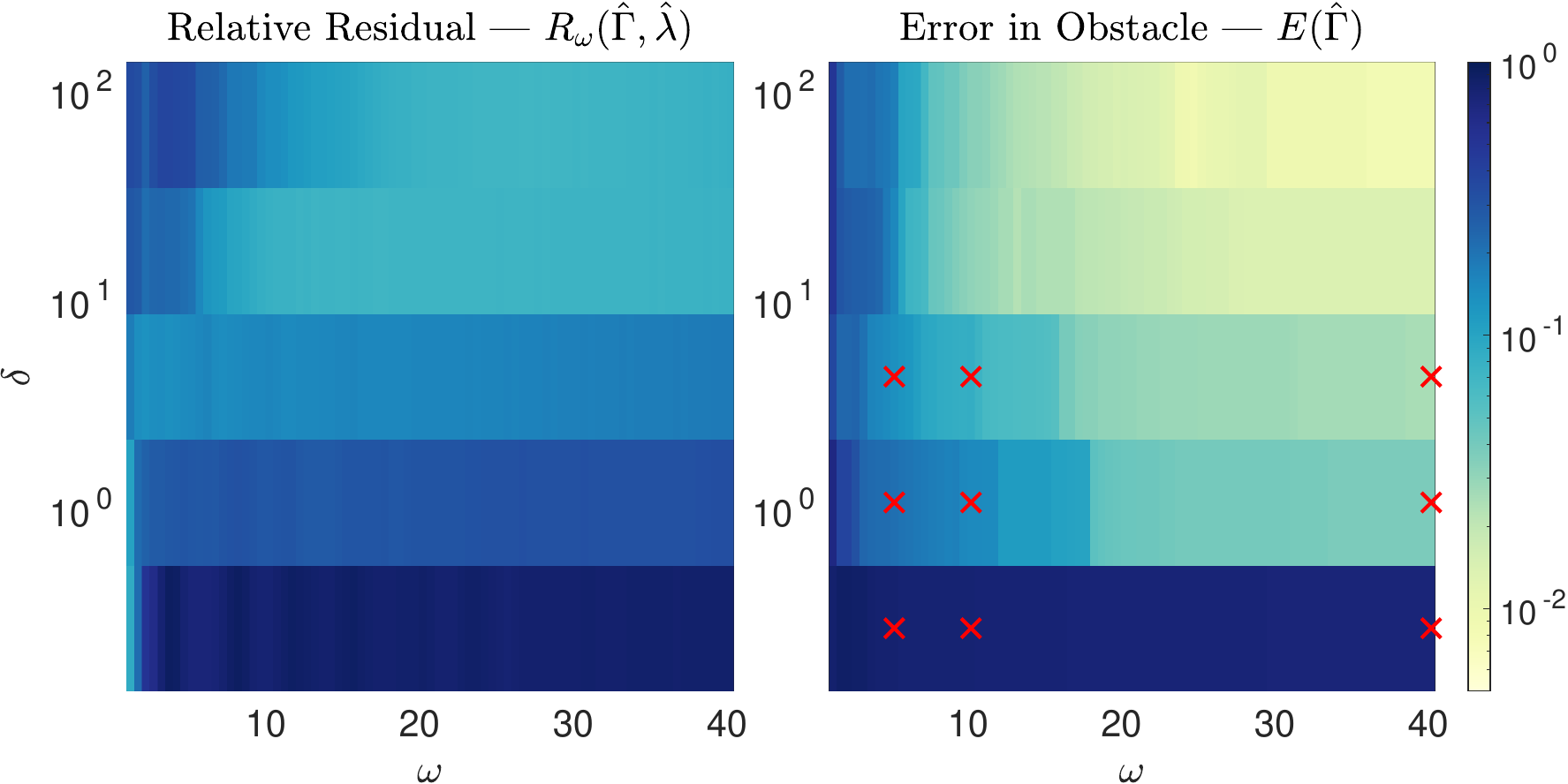}
  \caption{Experiment of \cref{sec:withnoise} with noise-corrupted transmission data
    using the $\lamabv$ impedance model for the inverse problem.
    The top row shows reconstructions obtained for different
    values of the pulsation, $\omega$, and dissipation, $\delta$. The bottom row has plots
    of the error measures as a function of $\omega$ and $\delta$.
    The values of $(\omega,\delta)$ for the reconstructions in the top
    row are marked in red in the error plot.}
    \label{fig:resanderr_varydeltafreq_noise1em1}
\end{figure}

Here we explore the effect of additive noise on the recovery of
the plane. We assume that the measured values are now of the
form

\begin{equation}
  \label{eq:noise_corrupted}
  \umeas_{\omega,\db} = \Ftrans_{\omega,\db}(\Gamma_\star) + \Sigma(\omega) \etab \; ,
\end{equation}
where the entries of $\etab \in \bbC^{N_r}$ are drawn from the standard
normal distribution and $\Sigma(\omega) = \sigma \max_\db |\Ftrans_{\omega,\db}(\Gamma_\star)|$
for some constant $\sigma$.

In this experiment, we generate transmission scattering data
for the more detailed plane with $\delta = \delta_0/4^j$ for
$j=0,1,2,3,4$. We then solve the inverse problem with the
$\lamabv$ model for each data set and $\sigma=10^{-1}$. The results are
shown in \cref{fig:resanderr_varydeltafreq_noise1em1}. 

The effect of the noise is visible in the
residual plot, where the smallest values of the residual are
similar to $\sigma$ and about an order of magnitude larger than in
the noise-free case (cf. \cref{fig:resanderr_varydeltafreq}).
On the other hand, the actual recovered domain performs at least
as well (and in some cases better) than in the noise-free case.
This suggests that the recovery is reasonably robust to additive
noise. 

\subsubsection{Experiment with backscatter data}

\label{sec:limitedaperture}

\begin{figure}
  \centering
  \includegraphics[width=0.5\textwidth]{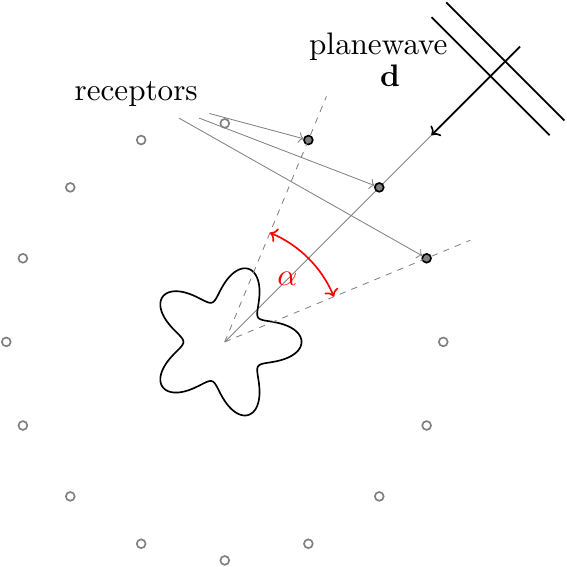}
  \caption{Illustration of arrangement of receptors for backscatter
    data.}
  \label{fig:backscatter-setup}
\end{figure}

\begin{figure}
  \centering
  \includegraphics[width=1\textwidth]{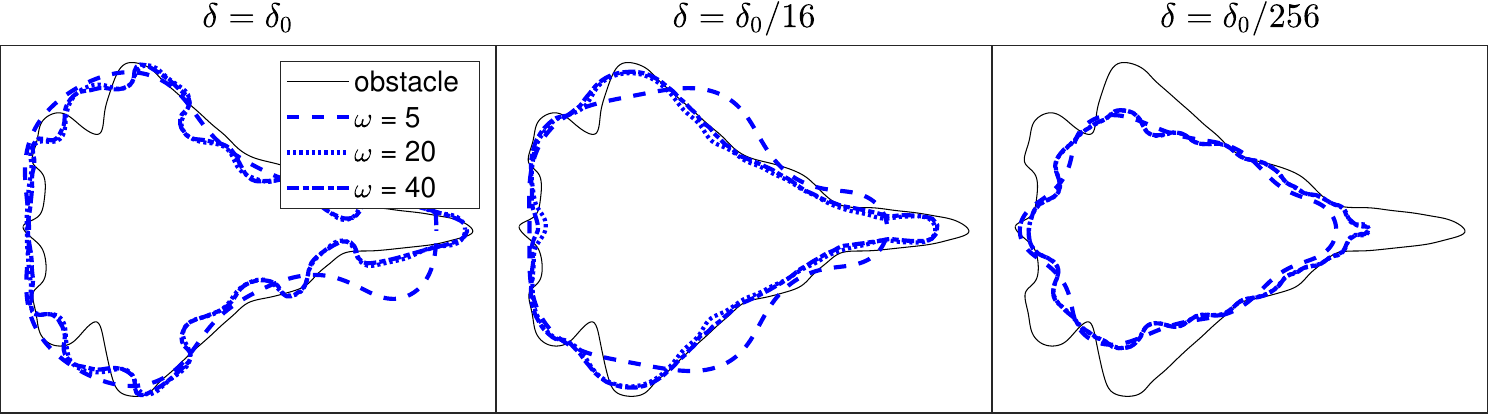}

  \vspace{0.5cm}
  
  \includegraphics[width=0.8\textwidth]{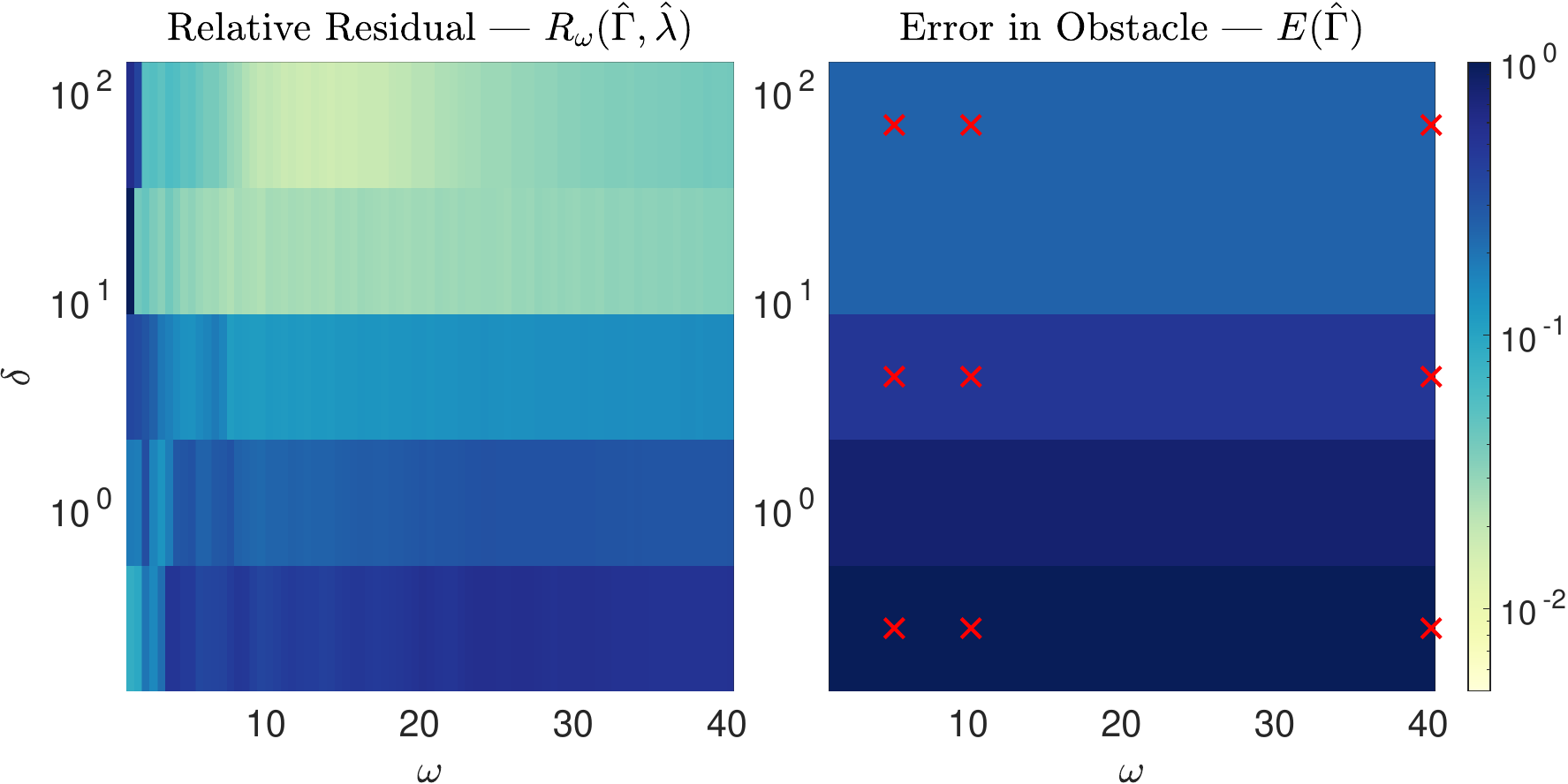}
  \caption{Experiment of \cref{sec:limitedaperture} with backscatter transmission data
    using the $\lamabv$ impedance model for the inverse problem.
    The top row shows reconstructions obtained for different
    values of the pulsation, $\omega$, and dissipation, $\delta$. The bottom row has plots
    of the error measures as a function of $\omega$ and $\delta$.
    The values of $(\omega,\delta)$ for the reconstructions in the top
    row are marked in red in the error plot.}
    \label{fig:resanderr_varydeltafreq_reflect}
\end{figure}

In the examples above, we have best-case scattering data, in that
the receptors surround the obstacle, there are many incident
directions, and the data are available at all receptors for
each incident wave. Such data may not be available in all
applications. Here we consider the recovery problem for backscatter
data. The backscatter set-up is illustrated in
\cref{fig:backscatter-setup}. For this problem, we have the same set
of possible receptor locations as for the experiments above, but
for each incident wave we assume that measurements are only
available within a certain angle ($\alpha$ in the figure) about
the axis defined by the planewave direction. For the experimental
results here, $\alpha = \pi/4$, so that the data set is $1/8$th the
size of the data set in \cref{sec:lowerdelta}.

The results for backscatter data are shown in \cref{fig:resanderr_varydeltafreq_reflect}.
The residual plot has minima at the lowest frequencies, where the
optimization procedure appears to get caught in local minima.
While the recovered domains bear some visual resemblance to
the true domain, the recovery is not particularly accurate for
any value of dissipation.

\subsubsection{Experiments with alternative models}

\label{sec:alternates}

\begin{figure}
  \centering
  \includegraphics[width=1\textwidth]{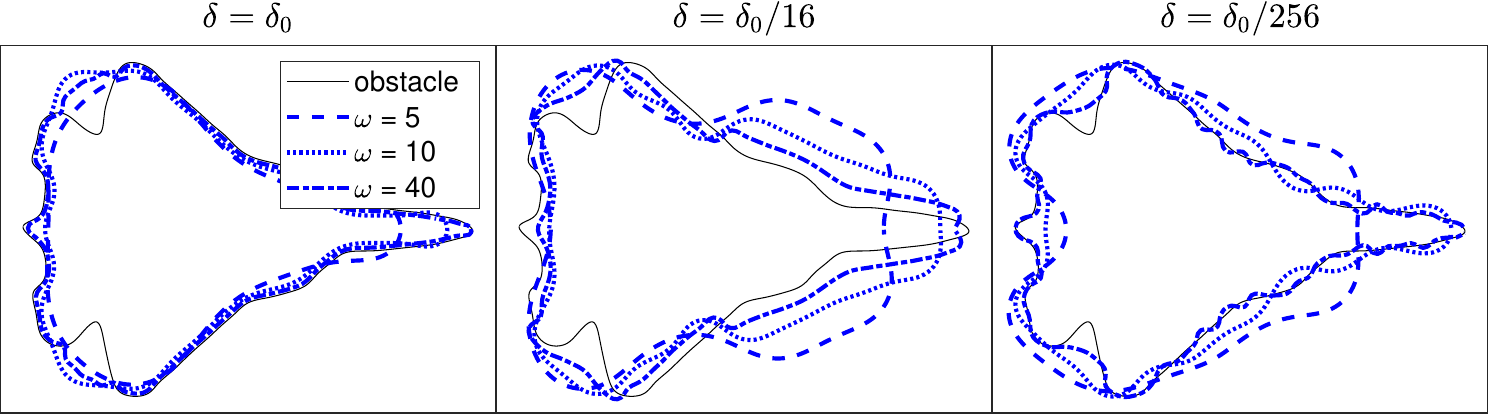}  
    
  \vspace{.5cm}
  
  \includegraphics[width=0.8\textwidth]{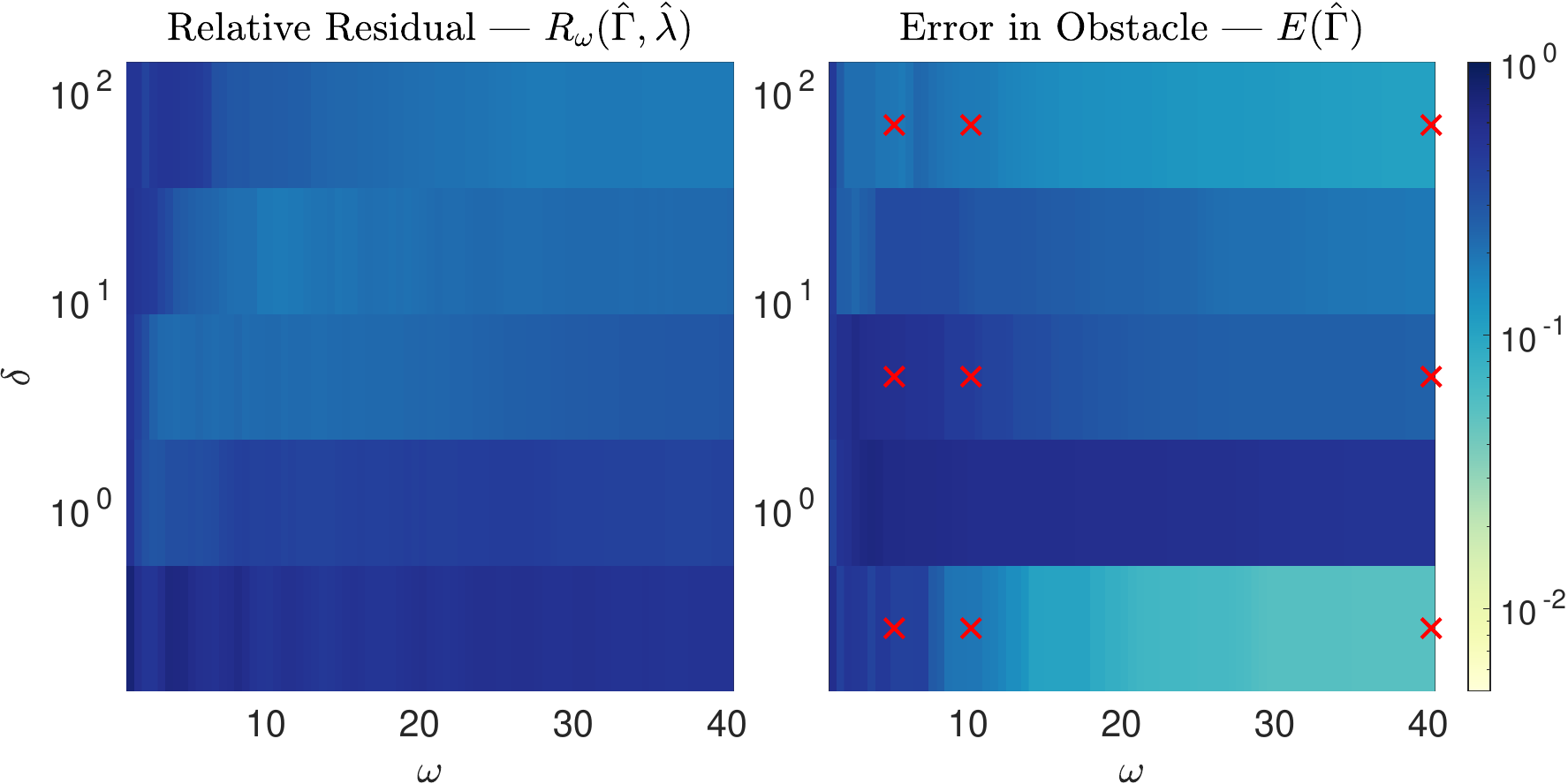}
  \caption{Experiment of \cref{sec:alternates} with transmission data
    using a constant impedance model for the inverse problem.
    The top row shows reconstructions obtained for different
    values of the pulsation, $\omega$, and dissipation, $\delta$. The bottom row has plots
    of the error measures as a function of $\omega$ and $\delta$.
    The values of $(\omega,\delta)$ for the reconstructions in the top
    row are marked in red in the error plot.}
    \label{fig:resanderr_varydeltafreq_constmodel}
\end{figure}

\begin{figure}
  \centering
  \includegraphics[width=1\textwidth]{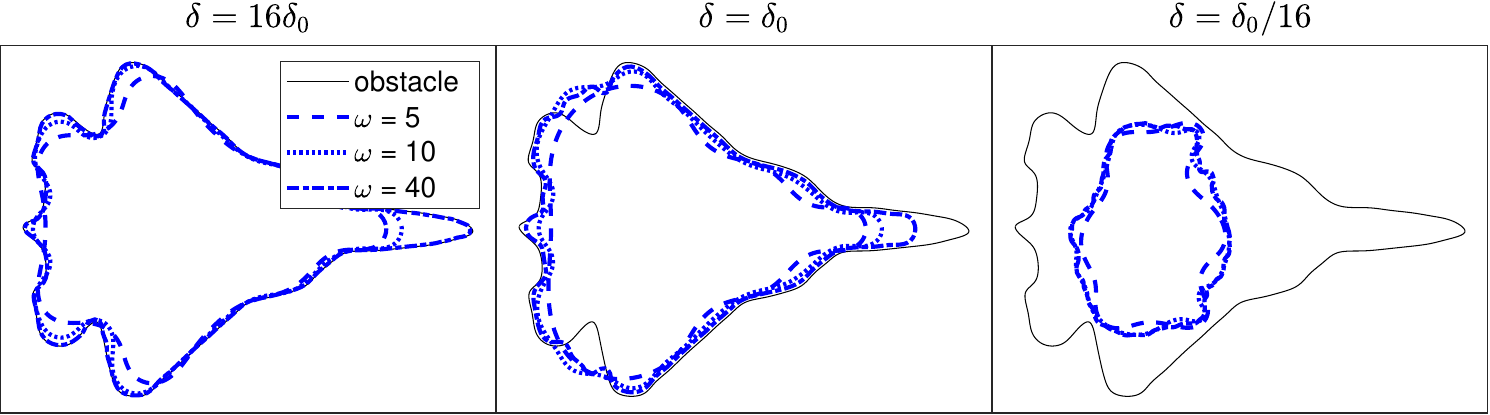}  
    
  \vspace{.5cm}
  
  \includegraphics[width=0.8\textwidth]{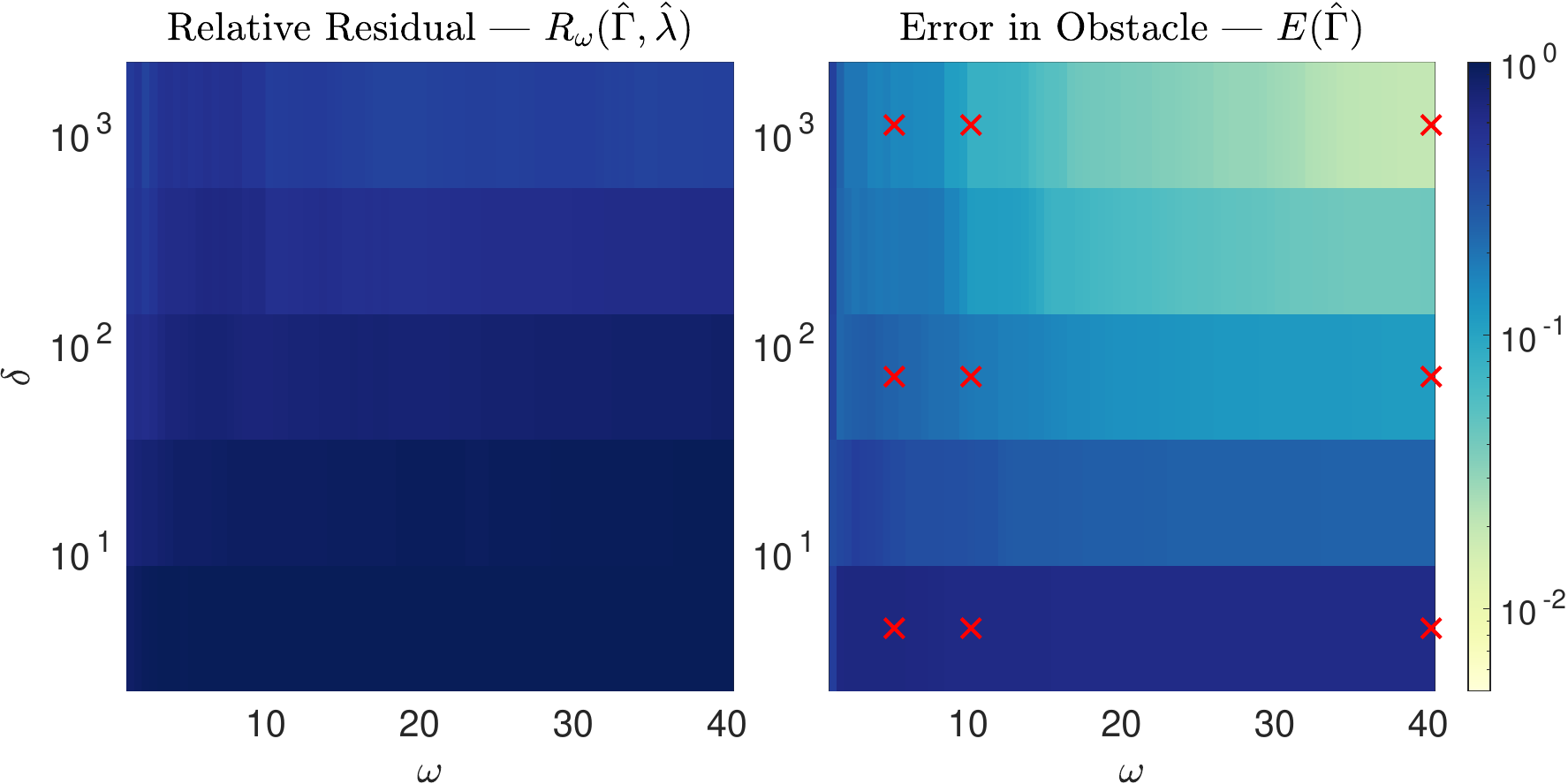}
  \caption{Experiment of \cref{sec:alternates} with transmission data
    using a sound-hard model for the inverse problem.
    The top row shows reconstructions obtained for different
    values of the pulsation, $\omega$, and dissipation, $\delta$. The bottom row has plots
    of the error measures as a function of $\omega$ and $\delta$. Note that
    the range of $\delta$ values is different from the other
    examples.
    The values of $(\omega,\delta)$ for the reconstructions in the top
    row are marked in red in the error plot.}
    \label{fig:resanderr_varydeltafreq_neumann}
\end{figure}

In these experiments, we verify that the curvature-dependent model
provides a meaningful advantage over even simpler models.

\paragraph{Constant impedance model}

Above, the failure mode of the curvature-dependent impedance models was addressed
by first seeking a minimum with the constant impedance model and then
using the obtained domain as an initial guess. In \cref{fig:const_model_initial_guess},
we saw that the constant impedance model resulted in a rather good
approximation of the domain. To explore the viability of this simpler model in
general, we provide the recovery results
for the constant impedance model applied to the same data as the experiment
of \cref{sec:lowerdelta} in \cref{fig:resanderr_varydeltafreq_constmodel}.

It is clear that the constant impedance model gets stuck in local minima
for even relatively large values of the dissipation. Compared with
\cref{fig:resanderr_varydeltafreq},
the constant impedance model performs worse than the $\lamabv$ model,
often by orders of magnitude. An interesting exception is the lowest
dissipation level considered ($\delta_0/256$), where the constant impedance
model results in a domain with some visual resemblance to the true domain but
the curvature-dependent model does not. 

\paragraph{Neumann boundary condition model}

As $\delta$ increases, we have that $\lamabv = \mathcal{O}(\delta^{-1/2})$. In the
limit, we then expect that the obstacle would behave as a sound-hard obstacle.
Here, we will test how well the sound-hard model does for the recovery of
the obstacle in the dissipative setting. 

Let $\Fneu_{\omega,\db}(\Gamma)$ be the vector of values of
$\phi(\rb_j)$ where $\phi$ is the solution of the PDE
\begin{equation}\label{eq:neupde}
\begin{aligned}
  -(\Delta + k_2^2)\phi &= 0 & \textrm{ in } \Omega_2 \; , \\ 
  \partial_n \phi &=
  - \partial_n \uinc \; & \textrm{ on } \Gamma \; , \\
  \sqrt{|\xb|} \left(\phi - ik_2\frac{\xb}{|\xb|}\cdot
  \nabla \phi \right) &\to 0  &\textrm{ as } |\xb|\to \infty \; .
\end{aligned}
\end{equation}
We will solve the minimization problem

\begin{equation}
  \label{eq:invprob_neu} 
 \hat{\Gamma}_j = \argmin_{\Gamma \in S_\Gamma(\omega_j)} \sum_{i=1}^{N_d} | \umeas_{\omega_j,\db_i} -
  \Fneu_{\omega_j,\db_i}(\Gamma)|^2 \; ,
\end{equation}
using similar tools to the ones we applied to the impedance
model. We consider a different range of values for the dissipation;
we set $\delta = \delta_0/4^j$ for $j=-2,-1,0,1,2$.

The recovery results are shown in \cref{fig:resanderr_varydeltafreq_neumann}.
While the Neumann model can recover the domain for the largest values of
the dissipation, the recovery obtained for $\delta \leq \delta_0$ has
significant error.

\section{Conclusion and future directions}\label{sec:conclusions}

We have developed an optimization framework for
using curvature-dependent impedance models, like $\lamabv$, to solve
inverse problems with transmission data from a dissipative obstacle.
The model succeeds in recovering the obstacle boundary reasonably
accurately, even at values of the dissipation for which waves transmit
through the narrowest parts of the obstacle. Intuitively, the model
fails for lower amounts of dissipation where many waves transmit.

The recovered material parameters are relatively 
accurate for the dissipation, $\delta$, and the product
of the relative sound speed and relative density, $c_r\rho_r$, but
the individual values of $c_r$ and $\rho_r$ appear more difficult to
recover. Because the domain is recovered well with the impedance
model, we plan to explore the use of a transmission model as a
post-processing step to obtain more accurate material parameters
once the domain is set. 

The model appears to be robust to (and may benefit from)
additive noise and provides a meaningful
advantage over simpler models, like constant impedance or Neumann
boundary conditions, in terms of the range of values of the dissipation
to which it applies.
Because an impedance model can emulate both
Neumann and Dirichlet boundary conditions in certain limits~\cite{borges2022multifrequency},
we believe that these results further support the use of a
curvature-dependent impedance model as a first-pass model for
any scattering data.

While the $\lamabv$ model was the focus of the numerical experiments,
the similar $\lamch$ model is more general and achieves results nearly as good. We plan to explore the performance of $\lamch$ in more settings;
for example, the more constrained $\lamabv$ is specific to the type of
transmission boundary condition considered here but $\lamch$ could
conceivably be used as a model for other transmission conditions and
thin coatings. 

Finally, the lackluster results with backscatter data suggest
that further work is needed for certain experimental settings.
Our experience is that manipulating the $M(\omega)$ and $\ccurv$
parameters was insufficient to improve the situation. We plan
to investigate improved regularization strategies and more
sophisticated optimization methods in these regimes.

\subsection*{Acknowlegments}The authors would like to thank Jeremy Hoskins and Manas Rachh for many useful discussions.

\subsection*{Funding}The work of C. Borges was supported in part by the Office of Naval Research under award number N00014-21-1-2389. 

\bibliography{references}

\appendix

\section{Integral equation formulations of the PDEs and their numerical discretization}
\label{sec:bieform}
To solve the forward impedance and transmission problems, we reformulate
the  partial differential equations, i.e.
\cref{eq:transpde} and
\cref{eq:imppde},
using well-established methods from layer potential theory~\cite{colton1998inverse}.
We briefly describe the representations and numerical methods
we use in this section. 

Consider an obstacle, $\Omega_1$, with a smooth boundary curve, $\Gamma$. 
Let $G^k(\xb,\yb)=\imath H^{(1)}_0(k\|\xb-\yb\|)/4$ be the Green's function for the
two dimensional Helmholtz equation and $\nb$ be the outward normal
to the obstacle boundary, $\Gamma$.
Let $S_k$ and $D_k$ denote the single and double layer operators,
respectively, i.e.
\begin{equation*}
S_k\sigma(\xb)=\int_{\Gamma} G^k(\xb,\yb) \sigma (\yb) ds(\yb), \quad 
D_k\sigma(\xb)=\int_\Gamma \frac{\partial G^k(\xb,\yb)}{\partial {\bf n}(\yb)}\sigma(\yb)ds(\yb) \; ,
\end{equation*}
for $\xb \in \bbR^2\setminus \Gamma$. For $\xb \in \Gamma$, we denote these
operators by $\cS_k$ and $\cD_k$, respectively. 
We denote the normal derivatives of $S_k$ and $D_k$ by
$\cK_k$ and $\cT_k$, respectively, i.e.
\begin{equation*}
\cK_k\phi(\xb) = \int_\Gamma \frac{\partial G^k(\xb,\yb)}{\partial {\bf n}(\xb)}\phi(\yb)ds(\yb), 
\end{equation*}
and
\begin{equation*}
\cT_k\phi(\xb) = \fp \int_\Gamma \frac{\partial^2 G^k(\xb,\yb)}{\partial {\bf n}(\xb)\nb(\yb)}\phi(\yb)ds(\yb) \; ,
\end{equation*}
where $\xb \in \Gamma$ and $\fp$ indicates that the integral is interpreted in the Hadamard
finite part sense.

On smooth curves, $\cS_k$, $\cD_k$, and $\cK_k$ are weakly singular
integral operators, while $\cT_k$ is hyper-singular. 

\subsection{Transmission problem}
\label{sec:bieform-trans}

To solve the transmission problem, we represent the scattered field inside the obstacle, denoted
$\uscat_{k_1}$, and outside the obstacle, denoted $\uscat_{k_2}$, as
\begin{equation}\label{eq:trans_combfield}
\begin{aligned}
  \uscat_{k_1} = -\uinc + \frac{1}{b_{1}} D_{k_{1}} \mu - \frac{1}{b_{1}} S_{k_{1}} \sigma \\
  \uscat_{k_2} = \frac{1}{b_{2}} D_{k_{2}} \mu - \frac{1}{b_{2}} S_{k_{2}} \sigma
\end{aligned}
\end{equation}
where $\mu$ and $\sigma$ are unknown densities defined on $\Gamma$ and
$b_1 = \alpha = 1/(\rho_r(1+\imath\delta/\omega))$ and $b_2 = 1$.

Applying the boundary conditions in \cref{eq:transpde} to the
representation \cref{eq:trans_combfield} we obtain the system of boundary integral
equations
\begin{equation}\label{eq:int_trans}
\begin{bmatrix} 
    I+\frac{1}{qb_2}\cD_{k_2}-\frac{1}{qb_1}\cD_{k_1} & -\left(\frac{1}{qb_2}\cS_{k_2}-\frac{1}{qb_1}\cS_{k_1}\right)  \\
    \cT_{k_2}-\cT_{k_1} & I-\left(\cK_{k_2}-\cK_{k_1}\right)    
\end{bmatrix} 
\begin{bmatrix} 
\mu\\
\sigma
\end{bmatrix} =
\begin{bmatrix} 
-\frac{1}{q} \uinc\\
-b_2\frac{\partial \uinc}{\partial n} 
\end{bmatrix} \; ,
\end{equation}
where $q = (1/b_1 + 1/b_2)/2$. This system is Fredholm and
has a unique solution for the wavenumbers treated in this paper.
While $\cT_{k_i}$ is hyper-singular for $i=1,2$, the operator
$\cT_{k_2}-\cT_{k_1}$ is weakly singular owing to cancellations in
the singularities. 

\subsection{Impedance problem}
To solve the impedance and Neumann ($\lambda = 0$)
problems, we use the regularized combined layer potential
representation proposed and analyzed in \cite{bruno2012regularized,turc2017well}. In particular,
we adopt the representation
\begin{equation}\label{eq:imp_combfield}
\uscat = \left(S_{k_2}+\imath{k_2}D_{k_2}\cS_{\imath|{k_2}|}\right) \phi \; .
\end{equation}

Applying the boundary conditions in \cref{eq:imppde} to the representation
\cref{eq:imp_combfield} and applying Calderon identities, gives the boundary integral
equation
\begin{multline}\label{eq:int_imp}
    \left[-\frac{2+\imath k_2}{4}I+\cK_{k_2}+
    \imath{k_2}\left(\left(\cT_{k_2}-\cT_{\imath|{k_2}|}\right)\cS_{\imath|{k_2}|}+\left(\cK_{\imath|{k_2}|}\right)^2\right)+\imath{k_2}\lambda\left(\vphantom{\frac{\imath{k_2}}{2}\cS_{\imath|{k_2}|}}S_{k_2}+ \right.\right. \\
    \left.\left. \imath{k_2}\cD_{k_2}\cS_{\imath|{k_2}|}+\frac{\imath{k_2}}{2}\cS_{\imath|{k_2}|}\right)\right]\phi = -\left(\frac{\partial \uinc}{\partial n}+\imath {k_2}\lambda\uinc \right)
    \; ,
\end{multline}
which is Fredholm and has a unique solution when $2+\imath{k_2}\neq0$.
Similar to the above, $\cT_{k_2}-\cT_{\imath|k_2|}$ is a weakly singular
operator, so that the calculations can be arranged with all integral kernels
used being weakly singular.

\subsection{Numerical Solution}
\label{sec:pdenumerics}
The Nystr\"om method can be used to solve both systems \eqref{eq:int_trans} and
\eqref{eq:int_imp}. We discretize the boundary $\Gamma$ using $N=\mathcal{O}(k)$,
with approximately 10 points per wavelength, where for the impedance boundary condition
$k=k_2$ and for the transmission problem $k=\max(k_1,k_2)$. To discretize the
weakly singular integral operators we apply the Hybrid Gauss-trapezoidal quadrature
rule of order 16 from
\cite{alpert1999hybrid}.

We invert the resulting linear system directly using Gaussian
elimination at the cost of $\mathcal{O}(k^3)$ operations. We store the inverse, so this
work is amortized over incident directions and the inverse can be applied in
$\mathcal{O}(k^2)$ operations per incident field to generate forward data.
With $\mathcal{O}(k)$ incident fields,
the total work at a given frequency is then $\mathcal{O}(k^3)$ for generating
data and for evaluating the objective function.
In the inverse problem, filling the Jacobian matrix, as described in
\cref{sec:opt_descent_directions},
requires $\mathcal{O}(k)$ PDE solves per incident direction resulting in
$\mathcal{O}(k^4)$ total work at a given frequency. These costs can be reduced by
using alternate optimization methods or by employing fast methods for the PDE
solutions but this was not a focus of the present work.

\end{document}